\documentclass[10pt,english,draftclsnofoot,onecolumn]{IEEEtran}
\IEEEoverridecommandlockouts

\usepackage{cite}
\usepackage{ifthen}
\usepackage{graphicx}
\usepackage{verbatim}
\usepackage{subcaption}
\usepackage{multirow}
\usepackage[usenames,dvipsnames]{color}
\usepackage{amssymb,amsmath,mathtools}
\usepackage{xcolor}
\usepackage{multirow}
\usepackage{url}
\usepackage{multicol}
\usepackage{amsmath}
\usepackage{dcolumn}
\usepackage{amsfonts}
\usepackage{latexsym}
\usepackage{graphicx,epstopdf}
\usepackage{color}
\usepackage{amsthm}
\usepackage{subeqnarray}
\usepackage{authblk}
\usepackage{indentfirst}
\usepackage{enumerate}
\usepackage{algorithm}
\usepackage{algpseudocode}
\makeatletter
\let\OldStatex\Statex
\renewcommand{\Statex}[1][3]{%
	\setlength\@tempdima{\algorithmicindent}%
	\OldStatex\hskip\dimexpr#1\@tempdima\relax}
\makeatother

\allowdisplaybreaks

\DeclareMathOperator{\N}{\mathcal{N}}

\newboolean{showcomments}
\setboolean{showcomments}{true}
\newcommand{\lina}[1]{  \ifthenelse{\boolean{showcomments}}
{ \textcolor{red}{(Lina says:  #1)}} {}  }
\newcommand{\Guannan}[1]{  \ifthenelse{\boolean{showcomments}}
{ \textcolor{blue}{(Guannan says:  #1)}} {}  }
\newcommand{\algoname}{\textproc{OptDist-VC}}
\newboolean{showhighlight}
\setboolean{showhighlight}{true}
\newcommand{\revise}[1]{\ifthenelse{\boolean{showhighlight}}{\textcolor{black}{#1} }{#1} }

\newboolean{isfullversion}
\setboolean{isfullversion}{true}

\def\ba{\begin{array}}
\def\ea{\end{array}}
\newcommand{\beq}{\begin{equation}}
\newcommand{\eeq}{\end{equation}}
\newcommand{\bq}{\begin{eqnarray}}
\newcommand{\eq}{\end{eqnarray}}
\newcommand{\bqn}{\begin{eqnarray*}}
\newcommand{\eqn}{\end{eqnarray*}}
\newcommand{\bee}{\begin{enumerate}}
\newcommand{\eee}{\end{enumerate}}
\newcommand{\bi}{\begin{itemize}}
\newcommand{\ei}{\end{itemize}}

\newcommand{\btab}{\begin{tabular}}
\newcommand{\etab}{\end{tabular}}

\newcommand{\qedd}{\ \hfill{$\Box$}\\[5pt]}
\newcommand{\R}{\mathbb{R}}
\newcommand{\LL}{\mathcal{L}}

\newcommand{\SSS}{\mathcal{S}}

\newcommand{\blambda}{\bar{\lambda}}
\newcommand{\ulambda}{\underline{\lambda}}

\newcommand{\bv}{\bar{v}}
\newcommand{\uv}{\underline{v}}
\newcommand{\barq}{\bar{q}}
\newcommand{\unq}{\underline{q}}

\newcommand{\tildex}{\tilde{X}}
\newcommand{\hatq}{\hat{q}}

\newcommand{\vv}[1]{\boldsymbol{\mathbf{#1}}}
\newtheorem{theorem}{Theorem}

\newtheorem{proposition}[theorem]{Proposition}
\newtheorem{lemma}[theorem]{Lemma}

\newtheorem{remark}{Remark}
\newtheorem{assumption}{Assumption}

\makeatother

\usepackage{babel}

\begin{document}

\title{Optimal Distributed Feedback Voltage Control under Limited Reactive Power}
\author{Guannan Qu, Na Li   \thanks{Guannan Qu and Na Li are with the School of Engineering and Applied Sciences, Harvard University, Cambridge, MA 02138, USA (Emails: gqu@g.harvard.edu, nali@seas.harvard.edu). This work was supported by NSF 1608509, NSF CAREER 1553407 and ARPA-E through the NODES program.
}}

\maketitle

\thispagestyle{plain}
\pagestyle{plain}

\begin{abstract}
	In this paper, we propose a distributed voltage control in power distribution networks through reactive power compensation. The proposed control can (i) operate in a distributed fashion where each bus makes its decision based on local voltage measurements and communication with neighboring buses, (ii) always satisfy the reactive power capacity constraint, (iii) drive the voltage magnitude into an acceptable range, and (iv) minimize an operational cost. We also perform various numerical case studies to demonstrate the effectiveness and robustness of the controller using the nonlinear power flow model. 
\end{abstract}
	\section*{{Nomenclature}}

\addcontentsline{toc}{section}{Nomenclature}

{ \color{black}
	\subsection{Parameters}
	
	\begin{IEEEdescription}[\IEEEusemathlabelsep\IEEEsetlabelwidth{$V_1,V_2,V_3$}]
		\item [$\mathcal{N},  n,\mathcal{E}$] $\mathcal{N}=\{0,1,\ldots,n\}$ is the set of buses with $0$ being the substation; $n$ is the number of buses (except the substation); $\mathcal{E}$ is the set of lines in the network.
		\item [$\mathcal{N}_i,\mathcal{P}_i$] Neighbors of bus $i$ in the network (including $i$ itself but excluding the substation); the set of lines in the network on the unique path from the substation to bus $i$.
		\item [$r_{ij}, x_{ij}$] Resistance and reactance on the transmission line between $i,j$.
\item [$X$, $R$] Matrix in linearized branch-flow model.
\item [$Y$]  Inverse of matrix $X$.
\item [$\vv{v}^{par}$]  Parameter in the linearized branch flow model that means the vector of voltage when no control action is taken.		
		\item [$\bar{q}_{i}, \underline{q}_{i}$] Upper, lower limits of reactive power injection at bus $i$.
		\item [$\bar{v}_{i}, \underline{v}_{i}, \vv{\bar{v}},\vv{\underline{v}}$] Upper, lower limits of voltage at bus $i$; vector of $\bar{v}_{i}$ and $\underline{v}_{i}$ at all buses.
		\item [$\mu,l$] Strongly convex and smooth parameter of the cost $f$.
		\item [$d$] Parameter in the cost function that balances the reactive power cost and the approximate power loss term.
		\item [$\alpha,\beta,\gamma,c$] Step size parameters in our controller.
	\end{IEEEdescription}	
	
	\subsection{Variables and Functions}
	\begin{IEEEdescription}[\IEEEusemathlabelsep\IEEEsetlabelwidth{$V_1,V_2,V_3$}]
		\item [$v_i,\vv{v}$] Squared voltage magnitude at bus $i$; vector of $v_i$ at all buses.
		\item [$q_i$, $\vv{q}$] Reactive power injection at bus $i$; vector of $q_i$ at all buses.
		\item [$p_i$, $\vv{p}$] Active power injection at bus $i$; vector of $p_i$ at all buses. 
		\item [$\hat{q}_i,\vv{\hat{q}}$] Virtual reactive power injection at bus $i$ or the primal variable in the Lagrangian; vector of $\hat{q}_i$ at all buses.
		\item [$v_i(\vv{q})$, $\vv{v}(\vv{q})$] Voltage at bus $i$ as a function of the reactive power injection across the buses; the vector of $v_i(\vv{q})$ at all buses. Not to be confused with $v_i(t), \vv{v}(t)$, the voltage measurement at time $t$.
		\item [$p_{ij},q_{ij},\ell_{ij}$]  Active power flow, reactive power flow and the square of current magnitude on the transmission line from $i$ to $j$.
		\item [$f_i, f$]  Cost of providing reactive power at bus $i$; the total cost. 
		\item [$\bar{\lambda}_i, \underline{\lambda}_i$] Lagrangian multiplier for the voltage upper limit and voltage lower limit at bus $i$. 
		\item [$\bar{\vv{\lambda}}$, $\underline{\vv{\lambda}}$, $\vv{\lambda}$]  Vector of $\bar{\lambda}_i$ at all buses; the vector of $\underline{\lambda}_i$ at all buses; the vector of $\bar{\lambda}_i$ and $\underline{\lambda}_i$ at all buses.
		\item [$\xi_i, \vv{\xi}$] Lagrangian multiplier for the reactive power capacity constraint at bus $i$; the vector of $\xi_i$ at all buses. 
		\item [$\mathcal{L}(\vv{\hat{q}}, \vv{\xi},\vv{\lambda})$]  Augmented Lagrangian for the optimization problem.  
	\end{IEEEdescription}
	\subsection{Notations}
	\begin{IEEEdescription}[\IEEEusemathlabelsep\IEEEsetlabelwidth{$\jmath:=\sqrt{-1}$}]
		\item[$\lbrack P\rbrack_{ij}$] The $i,j$-th entry of matrix $P$.
		\item[$\mathbf{1}$] $n\times1$ column vector with all ones. 	
		\item [$\R^{n}_{\geq 0}$] $n$-dimensional nonnegative orthant $\{ \vv{x} = [x_1,\ldots,x_n]^T: x_i\geq 0\}$.
		\item[$\lbrack \vv{x}\rbrack^+$ ] The projection of $\vv{x}$ onto the nonnegative orthant.
		\item [$\lbrack x\rbrack_{\underline{x}}^{\bar{x}}$] Projection of $x$ onto the interval $[\underline{x},\bar{x}]$, where $x$, $\bar{x}$, $\underline{x}$ are scalars.
		\item  [$\lbrack \vv{y}\rbrack_{\vv{\underline{y}}}^{\vv{\bar{y}}}$] Projection of $\vv{y}$ on to the box set $\{\vv{z}\in\R^n: z_i\in[\underline{y}_i,\bar{y}_i] \}$, where $\vv{y}$, $\vv{\underline{y}}$, $\vv{\bar{y}}$ are vectors of dimension $n$.
		\item [$\text{ST}_{b_1}^{b_2}(\cdot)$]  Soft thresholding function.
				\item[$\Vert \cdot \Vert$] Euclidean norm for vectors, spectral norm for matrices.
		\item[$\bar{\sigma}(P)$, $\underline{\sigma}(P)$] Largest, smallest eigenvalue of $P$. 
		\item [$\kappa(P)$] Condition number of matrix $P$. 

	\end{IEEEdescription}
	
}
\section{Introduction}
The primary purpose of voltage control is to maintain acceptable voltages at all buses along a distribution feeder under all possible operating conditions \cite{Baran1989a,Baran1989b,li2014real, zhang2013local, bolognani2013distributed, zhu2016fast,kekatos2015fast}.
Due to the increasing penetration of distributed energy resources (DER) such as photovoltaic and wind generation in the distribution networks, the operating conditions (supply, demand, voltages, etc) of the distribution feeder fluctuate fast and by a large amount. 
To overcome the challenges, it has been proposed to utilize the computing, (local) sensing, and (local) communication capabilities of the inverters in the DERs to adjust the reactive power injection in order to maintain voltage stability \cite{smith2011smart,turitsyn2011options}.


Various voltage control methods have been proposed \cite{farivar2012optimal,li2012exact,yeh2012adaptive,deshmukh2012voltage,farivar2013equilibrium,jahangiri2013distributed}.  One popular approach is using optimization methods.  Typically, an Optimal Power Flow (OPF) problem with voltage constraints is formulated and then ``computationally'' solved, either in a central (e.g. \cite{farivar2012optimal,li2012exact,yeh2012adaptive,deshmukh2012voltage}) or distributed way (e.g. \cite{zhang2015optimal,kraning2012message,dall2013distributed,vsulc2014optimal}). After obtaining the computational solution, the reactive power injection is adjusted.  
Hence, this approach is usually referred to as feedforward optimization: the algorithms assume knowledge of the disturbance (e.g. uncontrollable loads) and the system models while not using real-time measurements such as voltage magnitudes. In this paper, we mainly focus on distributed \textit{feedback} voltage control algorithms, in which the controller does not know the disturbance explicitly, but takes local measurements and adjusts its reactive power output based on the local measurements and local communication with its neighboring buses.  

There have been many efforts on developing feedback voltage control methods. One class is the traditional ``Droop'' control \cite{farivar2013equilibrium,jahangiri2013distributed}, as advocated by IEEE 1547-2018 Standard \cite{IEEE1547_2018}. It monitors the local bus voltages and adjusts the reactive power injection accordingly. However, \cite{li2014real} shows that droop type controllers are not able to maintain a feasible voltage profile under certain circumstances; ref. \cite[Section V-A]{zhu2016fast} shows droop control might experience stability and efficiency issues when the network is large. Therefore, other more sophisticated controllers have been proposed, e.g. \cite{li2014real, zhang2013local, bolognani2013distributed, zhu2016fast,kekatos2015fast, cavraro2016value,liu2018hybrid,tang2019fast}. Ref. \cite[Algorithm 1]{li2014real} and \cite{zhang2013local} propose an integration type controller that reaches a feasible voltage profile. Ref. \cite[Algorithm 2]{li2014real} and \cite{bolognani2013distributed} utilize a dual ascent approach that minimizes a power loss related cost while reaching a feasible voltage profile. However these methods also have their limitations. For instance, \cite{li2014real,zhang2013local,bolognani2013distributed} ignore the hard constraints on the reactive power injection capacity; \cite{zhu2016fast, kekatos2015fast} does not meet the hard voltage constraint.\footnote{Instead, \cite{zhu2016fast,kekatos2015fast} incorporate a weighted voltage deviation as a soft penalty.} Though \cite{cavraro2016value} meets both the voltage constraint and the hard reactive power constraint, there is a lack of theoretic guarantee for convergence even assuming linearized system model. 

Besides the concerns on convergence and voltage/reactive power constraints, another issue is the \textit{optimality} of voltage control. Since there is an acceptable range for voltage and reactive power, there is flexibility on the operating point of voltage control. Some operating points will have a lower operational cost than others. For example, for a DER with a fixed  apparent power rating, it is preferable for the DER to generate less reactive power so that it can generate more active power. In other cases, it may be preferable for DERs to operate at certain power factor, which requires its reactive power injection to be close to a certain value. Though some existing methods, e.g., \cite{li2014real, zhang2015optimal, bolognani2013distributed, cavraro2016value,zhu2016fast,kekatos2015fast}, do optimize a particular objective, the objective can not be freely chosen by DERs and does not necessarily reflect the true cost of DERs. It will be appealing if the voltage control method not only maintains the voltage in the acceptable range, but also minimizes a cost that reflects a meaningful operation cost. 


\textit{Our Contribution:} To overcome these challenges, we propose a distributed feedback voltage control that unifies the above controllers in the sense that it can simultaneously (i) meet the voltage constraint asymptotically, (ii) satisfy the reactive power capacity constraint throughout, and (iii) minimize an operation cost that can be composed of a power loss related cost and reactive power operation costs. The controller takes the voltage measurements as inputs and determines the reactive power injection through local communication and computation. The communication graph is the same as the physical distribution network, meaning that each bus only needs to communicate with its 1-hop neighbors.  The controller builds on the augmented Lagrangian multiplier theory \cite[Sec. 3.2]{bertsekas2014constrained} and primal-dual gradient algorithms \cite{nedic2009subgradient,cherukuri2016asymptotic,feijer2010stability,qu2018exponential}. We mathematically prove the performance of the controller using linear branch flow models \cite{Baran1989c} and numerically simulate the controller on a real distribution feeder using the {nonlinear power flow model}. We also test the {robustness} of the controller against measurement errors, communication delays, modeling errors and unbalanced networks. 
 
 We also note that the use of communication in our controller is inevitable. In fact, \cite{cavraro2016value} shows that for a class of communication-free controllers, there are scenarios in which those controllers can \textit{not} reach a feasible operating point (that satisfies both the voltage and the capacity constraint), despite the existence of a feasible operating point. The results in this paper are consistent with the performance limit in \cite{cavraro2016value} and demonstrate how to incorporate communication into controller design.
 
 {\color{black}
 	Finally, we note there are other related work \cite{massa2016dispersed,robbins2013two,carvalho2013distributed,kam2014stability,ali2012microgrid}. For example,   \cite{robbins2013two} proposes a online stochastic gradient descent based algorithm for online energy management; \cite{carvalho2013distributed} studies using active power curtailment to mitigate the voltage rise caused by DER. Also related is the large literature on microgrid control \cite{ali2012microgrid}. For example, \cite{massa2016dispersed} uses a centralized voltage control scheme to correct voltage deviation, where a large linear equation system coupling all buses needs to be solved.  
 }

\section{Preliminaries: Power Flow Model and Problem Formulation}\label{sec:problem}

\subsection{Branch flow model for radial networks and its linearization\label{sub:branchflow}}
In this paper, we consider a balanced distribution network with radial structure which consists of a set $\mathcal{N} = \{0,1,\dots,n\}$
of buses and a set $\mathcal{E}\in \mathcal{N}
\times \mathcal{N}$ of distribution lines connecting these buses.\footnote{The edges in $\mathcal{E}$ are directed, with the natural ordering that points towards the direction farther away from the substation. }
Bus $0$ represents
the substation and other buses in $\N$ represent branch buses. We let $\mathcal{N}_i$ denote the set of buses that are neighbors of bus $i$, including $i$ itself but excluding the substation bus. For
each line $(i,j)\in \mathcal{E}$,
 let $z_{ij}=r_{ij}+\textbf{i}x_{ij}$ be the impedance on line $(i,j)$, and $s_{ij}=p_{ij}+\textbf{i}q_{ij}$ be
the complex power flowing from buses $i$ to bus $j$. On each bus
$i\in \N$, let 
$s_{i}=p_i + \textbf{i} q_i$ be the
complex power injection. 
The branch flow model was first proposed in \cite{Baran1989a,Baran1989b} to model power flows in a radial distribution circuit \cite{Farivar-2012-BFM-TPS,Gan2012branch}: for each $(i,j)\in\mathcal{E}$,

{\small
\begin{subeqnarray}
\label{eq:branchflow_AC}
-p_j & = &  p_{ij} - r_{ij} \ell_{ij} - \sum_{k: (j,k)\in \mathcal{E}} p_{jk}
\label{eq:Kirchhoff.2a}
\\
-q_j & = &  q_{ij} - x_{ij} \ell_{ij}  - \sum_{k: (j,k)\in \mathcal{E}} q_{jk}
\label{eq:Kirchhoff.2b}
\\
v_j & = & v_i - 2 (r_{ij} p_{ij} + x_{ij} q_{ij}) + (r_{ij}^2 + x_{ij}^2) \ell_{ij}
\label{eq:Kirchhoff.2c}
\\
\ell_{ij} & = &  \frac{p_{ij}^2 + q_{ij}^2}{v_i}
\label{eq:Kirchhoff.2d}
\end{subeqnarray}}where $\ell_{ij}$ is the square of current magnitude on line $(i,j)$, and $v_i$ is the square of the voltage magnitude of bus $i$. As customary, we assume that the voltage magnitude $v_{0}$ on the substation bus is given and fixed at the nominal value (1 p.u.). 
Eq. (\ref{eq:branchflow_AC}) defines a system of equations in the variables $(p_{ij},q_{ij},\ell_{ij},v_i, (i,j)\in \mathcal{E})$.\footnote{Given these variables, the phase angles of voltages and currents can be uniquely determined for radial networks\cite{Farivar-2012-BFM-TPS,low2014convex}.}

{\color{black} There are many ways to obtain linearized power flow models, e.g. \cite{rueda2013mixed,garces2016linear,cespedes1990new}. In this paper, we use the Simplified Distflow model \cite{Baran1989c}, which, along with similar models, has been used in distributed voltage control in the literature \cite{zhu2016fast,liu2018hybrid,cavraro2016value,tang2019fast,bolognani2013distributed}. To obtain the Simplified Distflow model, we set the power loss term $\ell_{ij}$ to be $0$ in the branch flow model \eqref{eq:branchflow_AC} and get the following linear equations: for each $(i,j)\in\mathcal{E}$, }
\begin{subeqnarray}
\label{eq:branchflow_linear}
-p_j & = &  p_{ij} - \sum_{k: (j,k)\in E} p_{jk}
\\
-q_j & = &  q_{ij} - \sum_{k: (j,k)\in E} q_{jk}
\\
v_j & = & v_i - 2 (r_{ij} p_{ij} + x_{ij} q_{ij}) 
\end{subeqnarray} 

From (\ref{eq:branchflow_linear}), we can derive that the voltage vector $\vv{v}=[v_1,\ldots,v_n]^T$ and power injection vector $\vv{p}=[p_1,\ldots,p_n]^T,\vv{q}=[q_1,\ldots,q_n]^T$ satisfy the following equation: 
\begin{equation}
\vv{v}=R\vv{p}+X\vv{q}+v_{0}\vv{1}
\label{eq:v_pq}
\end{equation}
where $\vv{1}$ is a $n$-dimensional vector with all entries being $1$, and the $(i,j)$th entry of matrix $R,X\in\R^{n\times n}$ is given as follows 
{\small\begin{align*}
[R]_{ij} &:=2\sum_{(h,k)\in \mathcal{P}_i \cap \mathcal{P}_j} r_{hk}, \\
[X]_{ij}&:=2\sum_{(h,k)\in \mathcal{P}_i \cap \mathcal{P}_j} x_{hk}.  
\end{align*}}Here $\mathcal{P}_i \subset \mathcal{E}$ is the set of lines on the unique path from bus $0$ to bus $i$. The detailed derivation is given in \cite{lijun2013voltage}. In \cite{lijun2013voltage,bolognani2013distributed}, it has been shown that when the resistances and reactances of the lines in the network are all positive, the following proposition holds.
\begin{proposition}\label{prop:XY}
	If $x_{ij}>0$ for all edges $(i,j)$, then $X$ is positive definite; further define $Y:=X^{-1}$, then 
	\begin{align*}
	[Y]_{ij} = \left\{\begin{array}{ll}
	\sum_{k:k\sim i} x_{ik}^{-1}  & \text{ if } i=j\\
	-x_{ij}^{-1} & \text{ if } i\sim j \\
	0 & \text{ otherwise } 
	\end{array}  \right.
	\end{align*}
	here $i\sim j$ means $i\neq j$ and $i$ is a neighbor of $j$ in the network.
\end{proposition} 


\subsection{Problem formulation}\label{sec:prob}
We separate the reactive power $\vv{q}$ into two parts, $\vv{q} = \vv{q}^c+\vv{q}^{e} $, where $\vv{q}^c$ denotes the control action, i.e. the reactive power injection governed by the control components and $\vv{q}^e$ denotes any other reactive power injection. 
Let $\vv{v}^{par}:= R\vv{p}+X\vv{q}^e+v_0\vv{1}$. Then, given control action $\vv{q}^c$, the voltage profile by \eqref{eq:v_pq} is
\begin{align}
\vv{v}(\vv{q}^c)= X\vv{q}^c+\vv{v}^{par}. \label{eq:pf_inputoutput}
\end{align}
Here we intentionally write voltage profile $\vv{v}(\vv{q}^c)$ as a function of $\vv{q}^c$ to emphasize the input-output relationship between the control action and the voltage profile.
We comment that in \eqref{eq:pf_inputoutput}, $\vv{v}^{par}$ captures the loading conditions of the network that are \emph{unknown} and \emph{not controllable}. 
\revise{Without causing any confusion, we will simply use $\vv{q}$  instead of $\vv{q}^c$ to denote the reactive power injected by the control devices in the rest of paper.}

\textit{Feedback Control Loop.} We consider the following feedback control setting. 
At time $t$, let the reactive power injections of the nodes be $\vv{q}(t)$, then it determines the voltage profile $\vv{v}(t)\coloneqq \vv{v}(\vv{q}(t))$ through (\ref{eq:pf_inputoutput}).\footnote{Without causing any confusion, we abuse the notation $\vv{v}(\cdot)$ to denote both the network model (\ref{eq:pf_inputoutput}) mapping $\vv{q}$ to the voltage profile $\vv{v}(\vv{q})$, and the voltage profile at a certain time step $\vv{v}(t) $. }  Then given the voltage profile $\vv{v}(t)$ and other available information, the controller determines a new reactive power injection $\vv{q}(t+1)$. Mathematically, the voltage control problem is formulated as the following closed loop dynamical system,
\begin{align*} 
\vv{v}(t) & =  \vv{v}(\vv{q}(t)) = X\vv{q}(t)+\vv{v}^{par} \\ 
q_i(t+1) & =  \text{Controller}_i(\text{information available to $i$ at time } t) 
\end{align*}
{\color{black}where the information available to $i$ at $t$ may include the local voltage measurement, information received from neighbors, and locally stored internal states (that may have their own update equations).  }

\begin{figure}
	\centering
	\includegraphics[scale=0.4]{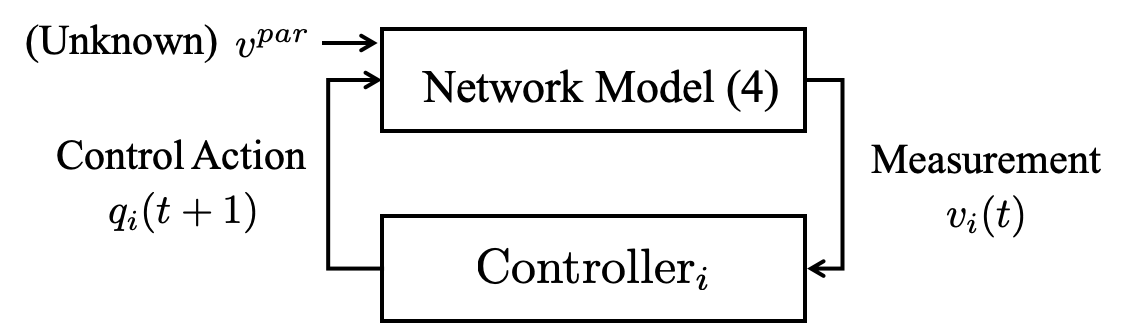}
	\caption{Feedback structure of the voltage control problem.}\label{fig:feedback}
\end{figure}
We also illustrate the feedback structure of the problem in Figure~\ref{fig:feedback}. We make the following remarks regarding the feedback nature of our problem.
\begin{itemize}
	\item The network is treated as an input-output system, and the controller has no knowledge of the internal details of the model (e.g. the $\vv{v}^{par}$ in \eqref{eq:pf_inputoutput}), except that it can measure the output $\vv{v}$ of the model. For now we assume the model is the Linearized DistFlow \eqref{eq:pf_inputoutput}, because the particular structure of \eqref{eq:pf_inputoutput} will facilitate our controller design and theoretic performance analysis. We emphasize that the controller in this paper can also be applied to  nonlinear power flow model \eqref{eq:branchflow_AC}. In our numerical case studies, we test the performance of our controller on the \textit{nonlinear} model.
	\item To facilitate theoretic analysis, we assume the $\vv{v}^{par}$(which captures the uncontrollable active and reactive generations and demands) in \eqref{eq:pf_inputoutput} is fixed during the control process. This is not required when applying our feedback controller. In the numerical case study, we test our controller under time varying loading conditions to show how our controller automatically adapt to the changing conditions. 
\end{itemize}



 \textit{Voltage Control Problem.} 
The objective of voltage control is to design a controller that meets the following four requirements.

\vspace{6pt}
\noindent\textbf{Requirement 1: Information.} The controller at $i$ shall only use information that is accessible locally and from neighboring buses in the network. This includes local decision variable $q_i(t)$, local voltage measurement $v_i(t)$, other local auxiliary variables and variables communicated from neighboring buses. 

\vspace{6pt}
\noindent\textbf{Requirement 2: Asymptotic voltage constraint.} \revise{The voltage profile reaches the acceptable lower and upper limits $\uv_i,\bv_i$, i.e. $\forall i$, $v_i(t)$ converges to a point inside the interval $[\uv_i,\bv_i]$.}

\vspace{6pt}
\noindent\textbf{Requirement 3: Hard capacity constraint.}  For each $i$, we introduce scalar $\underline{q}_i$ and $\bar{q}_i$, the lower and upper reactive power capacity limit for the device at node $i$. We require that, for any $t$, this capacity constraint shall \textit{not} be violated, i.e. $\unq_i\leq q_i(t)\leq \barq_i$. 

\vspace{6pt}
\noindent\textbf{Requirement 4: Optimality.}  We introduce $f_i :\mathbb{R}\rightarrow\mathbb{R}$, the operating cost of control action for each individual bus $i$.
We require that, under any system condition $\vv{v}^{par}$, the controller drives the distribution system to the optimal point of the following optimization problem, 
{\small\begin{subeqnarray}\label{eq:opt}
\min_{q_i} & & f(\vv{q}) \triangleq \sum_{i=1}^n f_i (q_i) + \frac{d }{2} \vv{q}^T X \vv{q} \slabel{eq:opt_cost}\\
s.t.  &  &\underline{v}_i \leq v_i(\vv{q}) \leq \bar{v}_i \slabel{eq:vol_con}\\
& &\underline{q}_i \leq q_i \leq \bar{q}_i \slabel{eq:q_con}
\end{subeqnarray}}In the optimization problem, the cost function (\ref{eq:opt_cost}) is composed of the sum of the operating costs $f_i$, as well as a network level cost $\frac{1}{2} q^T X q$, with a weighting parameter $d\geq 0$ balancing the two costs. {\color{black}The cost $\frac{1}{2} q^T X q$ is an approximation network loss term (up to a multiplicative factor and an additive term that does not depend on $\vv{q}$), under the assumption that the $R/X$ ratio of the network is constant \cite[Lemma 2]{bolognani2013distributed}. This cost has appeared in a few related work \cite{cavraro2016value,li2014real}. If the cost $\frac{d}{2} q^T X q$ is not an important factor for the system operator, then $d$ can be set as $0$.}
 We make the following regularity assumption on problem~(\ref{eq:opt}).
\begin{assumption}\label{assump:opt}
	(i) The cost function $f$ is differentiable, and is $\mu$-strongly convex and $l$-smooth, i.e. $\forall \vv{q},\vv{q}'\in\R^n$,  
	$$ \mu \Vert \vv{q}-\vv{q}'\Vert^2\leq\langle \vv{q} - \vv{q}', \nabla f(\vv{q}) - \nabla f(\vv{q}')\rangle\leq l \Vert \vv{q}-\vv{q}'\Vert^2. $$ 
	(ii) There exists a feasible solution $\vv{\hatq}^0$ for problem~\eqref{eq:opt} that meets the voltage constraint \eqref{eq:vol_con} with strict inequality. In other words, $\vv{\hatq}^0$ satisfies $\underline{v}_i < v_i(\vv{\hatq}^0) < \bar{v}_i$ and $\underline{q}_i \leq \hatq_i^0 \leq \bar{q}_i$, $\forall i$.
\end{assumption}
\begin{remark}  
Assumption~\ref{assump:opt}(i) is a standard assumption in the optimization literature. For the particular cost function in \eqref{eq:opt_cost} to meet Assumption \ref{assump:opt}(i), we need $f_i$ to be convex with Lipschitz gradients. For example, any linear or quadratic $f_i$, or zero cost function $f_i =0$ meets Assumption~\ref{assump:opt}(i). Assumption~\ref{assump:opt}(ii) ensures that the problem we are considering is feasible. 
\end{remark}
{\color{black}
	\begin{remark}\label{rem:cost}
		The reactive power provider sets the cost function $f_i$. One example of $f_i$ is the loss of opportunity cost for providing active power service. If the apparent power limit for the device is $s_i^{max}$, and the provider can sell the active power at price $\eta$, then the provider can earn (at most) $\eta \sqrt{(s_i^{max})^2 - q_i^2}$ when generating $q_i$ amount of reactive power. If the provider does not generate any reactive power, the provider can earn $\eta s_i^{max}$.  Therefore, the lost of opportunity cost is \begin{align*}
		\eta s_i^{max} - \eta \sqrt{(s_i^{max})^2 - q_i^2} &= \frac{\eta^2 q_i^2}{\eta s_i^{max} +\eta \sqrt{(s_i^{max})^2 - q_i^2}} \\
		&\approx \frac{\eta}{s_i^{max}} q_i^2 
		\end{align*} 
		which can be set as the cost function $f_i(q_i)$. Intuitively, the above cost function can also be understood as a penalization on large reactive power injection in order to provide more flexibility in supplying active power or maintain a power factor closer to $1$. 
		 Finally, we would like to emphasize that the fact our algorithm can guarantee optimality is an add-on functionality on top of the more basic requirements of meeting the voltage constraint with limited reactive power and local communication. If the users have no cost in supplying $q_i$, $f_i$ can be set to $0$.
	\end{remark}

}
{\color{black}
\begin{remark}
For easy exposition and without loss of generality, we assume there is a control component at each bus $i$. 
In Section~\ref{subsec:subcontrollable}, we discuss the scenario when some buses do not have control components.
\end{remark}}



\section{Distributed Voltage Controller}\label{sec:alg}
 In this section, we formally introduce our controller, Optimal Distributed Feedback Voltage Control (\algoname). For each bus $i$, we introduce auxiliary variables, $\hat{q}_i, \xi_i,\blambda_i,\ulambda_i$.  At each iteration $t$, node $i$ measures the local voltage $v_i(t)$,  then computes variables $\hatq_i(t+1)$,  $q_i(t+1)$, $\xi_i(t+1)$, $\blambda_i(t+1)$, $\ulambda_i(t+1)$ and injects the reactive power $q_i(t+1)$, and lastly passes certain variables to its neigboring buses in the network. The detailed implementation of \algoname\ is given as follows.
 
 \bigskip
\noindent \algoname:

\noindent {At time $t$, each bus $i$ follows the following $4$ steps.}

\noindent	\textbf{Step 1 (Measuring): } Receive voltage measurement $v_i(t)$.

\smallskip
\noindent	\textbf{Step 2 (Calculating):} Calculate $\hatq_i(t+1),\xi_i(t+1),\blambda_i(t+1),\ulambda_i(t+1)$ as follows.

\begin{figure}[t]
	\begin{center}
			\includegraphics[scale=0.5]{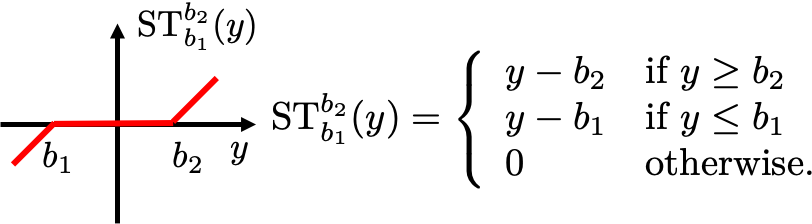}
	\end{center}
	\caption{The soft thresholding function.}\label{fig:st}
\end{figure}

{\small  \begin{subeqnarray}\label{eq:algorithm}
 	\hat{q}_i(t+1) & = & \hat{q}_i(t) - \alpha  \bigg\{ \bar{\lambda}_i(t)  - \underline{\lambda}_i(t) + d \hatq_i(t)  \nonumber \\
 	&&  + \sum_{j\in \mathcal{ N}_i} [Y]_{ij} \Big[f_j'(\hat{q}_j(t )  ) + \text{ST}_{c \underline{q}_j}^{c\bar{q}_j}(\xi_j(t)+ c \hatq_j(t)) \Big]  \bigg\}\nonumber \\
 	&& \slabel{eq:algorithm_hatq} \\
 	\xi_i(t+1) &=& \xi_i(t) + \beta \frac{\text{ST}_{c \underline{q}_i}^{c\bar{q}_i}(\xi_i(t)+ c \hatq_i(t))   - \xi_i}{c}\slabel{eq:algorithm_xi} \\
 	\blambda_i(t+1) &=& [\blambda_i(t) + \gamma(v_i(t) - \bar{v}_i) ]^+ \slabel{eq:algorithm_blambda} \\
 	 	\ulambda_i(t+1) &=& [\ulambda_i(t) + \gamma( \underline{v}_i - v_i(t) ) ]^+ \slabel{eq:algorithm_ulambda}
 	\end{subeqnarray}}where $[\cdot]^+ $ means projection onto the nonnegative orthant; quantity $\alpha,\beta$, $\gamma$ and $c$ are positive scalar parameters. For any $b_1<b_2$, function $\text{ST}_{b_1}^{b_2}(\cdot)$ is the soft-thresholding function defined as, $\text{ST}_{b_1}^{b_2}(y) = \max(\min(y - b_1,0),y-b_2)$ (see Fig.~\ref{fig:st} for an illustration).
  
  \smallskip
  \noindent\textbf{Step 3 (Injecting Reactive Power): }
Set reactive power injection at time $t+1$ as \begin{align}
  q_i(t+1) = [\hat{q}_i(t+1)]_{\unq_i}^{\barq_i}  \label{eq:algorithm_q}
 \end{align} 
 where  $[\cdot]_{\unq_i}^{\barq_i} $ means projection onto the set $[\unq_i,\barq_i]$. 
 
 \smallskip
\noindent\textbf{Step 4 (Communicating):}  	 Send values $f_i'(\hat{q}_i(t+1 )  ) + \text{ST}_{c \underline{q}_i}^{c\bar{q}_i}(\xi_i(t+1)+ c \hatq_i(t+1))$ to neighbors.\qedd

\begin{figure}
	\begin{center}
	\includegraphics[width=0.7\columnwidth]{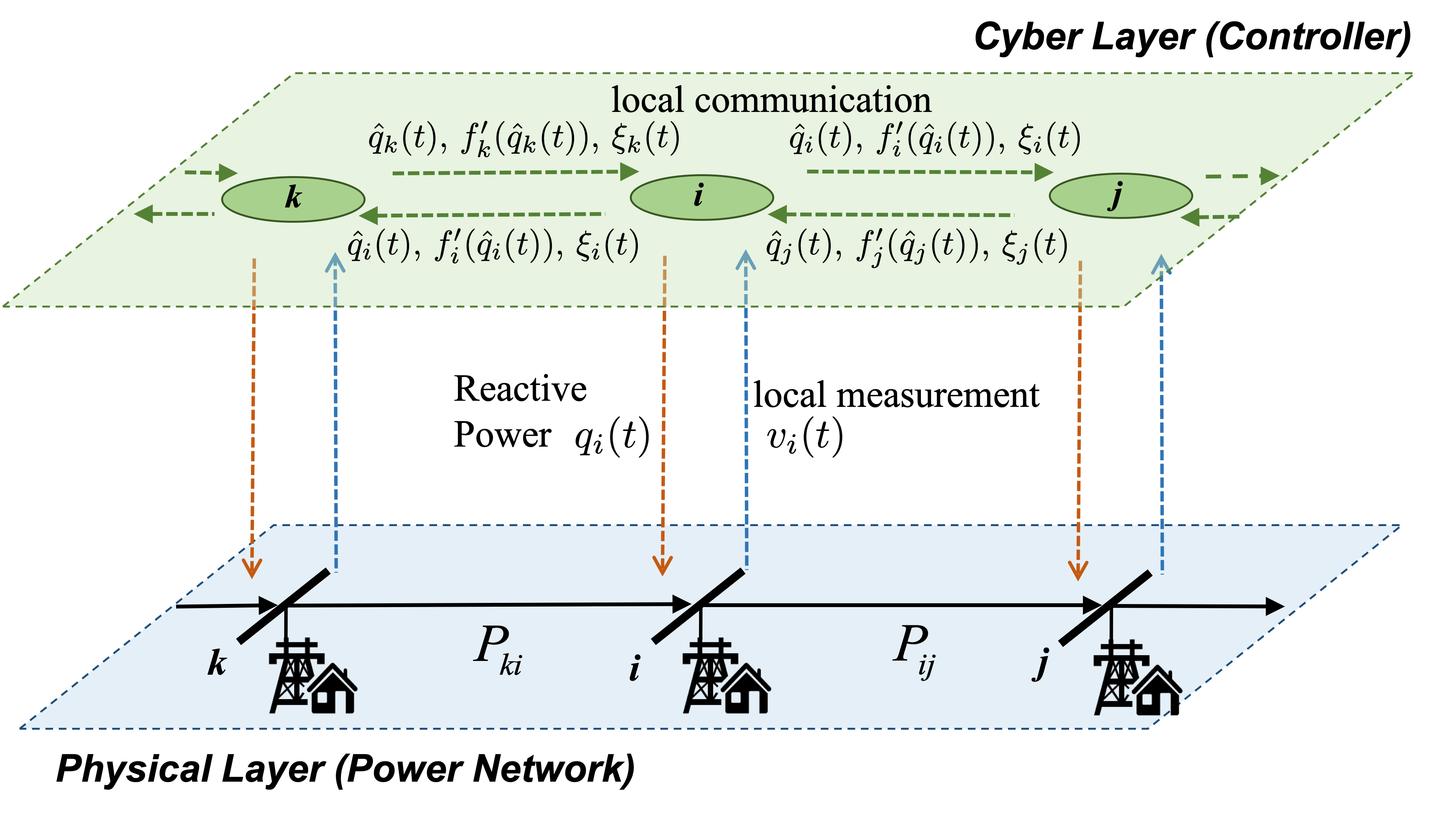}
	\end{center}
	\caption{Information Flow of \algoname.}\label{fig:algo}
\end{figure}

Figure~\ref{fig:algo} shows the information exchange between different buses and between the cyber layer (controller) and physical layer (network model) under \algoname. As Figure~\ref{fig:algo} shows, the cyber layer interacts with the physical layer only through voltage measurement $v_i(t)$ (Step 1) and reactive power injection $q_i(t)$ (Step 3), while all other activities of \algoname\ are conducted entirely inside the cyber layer, including calculation (Step 2) and communication (Step 4). A few comments on \algoname\  are in place. 
	\begin{itemize}
		\item $q_i(t)$ and $v_i(t)$ are physical quantities (reactive power injection and voltage), while ($\hatq_i(t), \xi_i(t),\blambda_i(t),\ulambda_i(t)$) are ``digital'' variables stored in the memory of the controller.
		\item Variable $\hatq_i(t)$ is the desired amount of reactive power to inject at time $t$. However it might violate the capacity constraint. Therefore, we simply set $q_i(t)$ as the projection of $\hatq_i(t)$ onto the capacity constraint. As will be shown later, $\hatq_i(t)$ will eventually become very close to $q_i(t)$.
		\item The update equation (\ref{eq:algorithm_hatq}) for the desired reactive power injection $\hatq_i(t)$ is to drive $\hatq_i(t)$ towards the superposition of the gradient of $f$ and certain ``correction directions'', related to $\blambda_j(t)-\ulambda_j(t)$ and $\xi_j(t)$, which pulls $\hatq_i(t)$ into the constraints set. Because of the superposition of the two directions, $\hatq_i(t)$ will be driven to minimize $f$ and in the meanwhile avoid violating the constraints.
		\item The variable $\xi_i(t)$, $\blambda_i(t),\ulambda_i(t)$ are Lagrangian multipliers that reflect the level of constraint violation of $\hatq_i(t)$. 
	\end{itemize}


It can be readily seen that \algoname\ meets Requirement 1 since it only needs local and neighboring information. Moreover, (\ref{eq:algorithm_q}) implies that $q_i(t)$ always lies within $[\unq_i,\barq_i]$ and hence Requirement 3 is met. At last, the following Theorem~\ref{thm:convergence} shows that $\vv{q}(t)$ will converge to the optimal solution of (\ref{eq:opt}) and hence Requirement 2 and 4 are satisfied. \textit{In conclusion, \algoname\ meets all the four design requirements.}

\begin{theorem}\label{thm:convergence}
	In \algoname, for any $c>0$, when $\alpha$, $\beta$, $\gamma$ are small enough and satisfy mild conditions,\footnote{For the detailed conditions on the step sizes, please see \ifthenelse{\boolean{isfullversion}}{Remark~\ref{rem:stepsize} in Appendix-\ref{appendix:proof_convergence}.}{Remark 6 in \cite[Appendix-B]{fullversion}.}} $\vv{q}(t)$ will converge to the unique optimizer of (\ref{eq:opt}).
\end{theorem}

 \algoname\ is based on a primal-dual gradient algorithm \cite{nedic2009subgradient,cherukuri2016asymptotic,feijer2010stability,qu2018exponential} for an augmented Lagrangian \cite{bertsekas2014constrained}, in which $\hat{q}_i(t)$ is the primal variable, $\xi_i(t),\blambda_i(t),\ulambda_i(t)$ are the dual variables. 
 We will elaborate on this in Section~\ref{subsec:concept}.
 
 {\color{black}We here provide a comparison of \algoname\ with related methods in \cite{li2014real, zhang2013local, bolognani2013distributed, zhu2016fast, cavraro2016value,liu2018hybrid,tang2019fast}. In terms of underlying methodology, \cite[Algorithm 1]{li2014real}, \cite{zhu2016fast} are based on (projected) primal gradient descent, \cite[Algorithm 2]{li2014real}, \cite{bolognani2013distributed,cavraro2016value,tang2019fast} are based on dual gradient ascent, while \cite{liu2018hybrid} and our method are based on (partial) primal-dual gradient algorithm. Compared to our method, \cite{li2014real,bolognani2013distributed} do not consider reactive power capacity limit; \cite{zhu2016fast,liu2018hybrid,tang2019fast} treat the voltage constraint as a soft penalty in the cost function in the form of (weighted) voltage deviation, whereas our work treats the voltage constraint as hard constraint; \cite{cavraro2016value} does not have theoretic guarantee for convergence when the reactive power constraint is enforced at every time step.\footnote{Ref. \cite{cavraro2016value}
 	does have theoretic guarantee for a variant of the algorithm therein, but that variant may violate the capacity constraint before it converges to an equilibrium point. } We provide a summary of the comparison of these methods in Table~\ref{tab:comparison}. }
 \begin{table*}
 	\centering
{\color{black} 	\caption{Comparison with related methods. } \label{tab:comparison}
 	\begin{tabular}{c|ccccc}
 		Algorithm & Underlying Methodology  &Voltage Constraint & Reactive Power Constraint & Cost & Communication  \\ \hline
 		\cite[Algo. 1]{li2014real} & Integral control & Satisfied  & Not considered  & N/A& Not needed   \\
 		\cite[Algo. 2]{li2014real} & Dual ascent & Satisfied & Not considered & $\frac{1}{2}q^T X q$ & Not needed   \\
 		\cite{zhang2013local} & Integral control & Penalization in cost  & Not considered  & $v$ deviation& Not needed   \\
 		\cite{bolognani2013distributed} & Dual ascent & Satisfied & Not considered & Power loss & Local comm.  \\
 		\cite{zhu2016fast} & Gradient Descent &  Penalization in cost  & Considered & Weighted $v$ deviation + $q$ cost  & Not needed \\
 		\cite{cavraro2016value} & Dual ascent & Satisfied & Considered & $\frac{1}{2}q^T X q$ & Local comm. \\
 		\cite{liu2018hybrid}& Primal-Dual & Penalization in cost & Considered & Weighted $v$ deviation + $v$ deviation & Local comm. \\
 		\cite{tang2019fast} & Accelerated dual ascent& Penalization in cost & Considered & $v$ deviation & Local comm. \\
 		Our work & Primal-Dual & Satisfied & Considered & $\frac{1}{2}q^T X q + q$ cost & Local comm.\\
 	\end{tabular}}
 	
 \end{table*} 
  {\color{black}
  	\begin{remark} The values of $\alpha,\beta,\gamma$ will affect convergence, and hence controller performance. One one hand, these values should be relatively small to make the algorithm converge. On the other hand, very small $\alpha,\beta,\gamma$ will make the convergence very slow.  We provide a theoretic upper bound on $\alpha,\beta,\gamma$ in Remark~\ref{rem:stepsize}. However, the theoretic upper bound is conservative, which is due to all known step size bounds that guarantee the ``contraction'' property of the primal-dual gradient algorithm, a property that is crucial in the analysis, are known to be conservative. For this reason, in the simulations we use trial and error to obtain the step sizes. Roughly speaking, the step size values would depend on the network structure and parameters, as well as the condition number of $f$.  In practice, to choose good parameter values would require some modeling and numerical testing before physical implementation, which is common for many (control) algorithms. {\color{black} Finally, we comment that our current results assume the network topology is fixed. It remains our future work to address the issue of robustness against system reconfiguration.}
  	\end{remark}}
  {\color{black}
  \begin{remark}
  	Under the IEEE 1547-2018 standard \cite{IEEE1547_2018}, our algorithm fits into the voltage-reactive power operation mode, where the reactive power injection responds to the voltage measurement. 
  \end{remark}
  	 }

 \section{Rationale Behind Controller Design}\label{subsec:concept}
  In this section, we will describe in detail the rationale behind our controller design, whereas the detailed mathematical proof \ifthenelse{\boolean{isfullversion}}{is deferred to the appendix.}{is deferred to the appendix of our online report \cite{fullversion} due to the space limit.}
 \subsection{\algoname\ as a primal-dual algorithm} \label{subsec:algo_derivation}
We introduce Lagragian multipliers $\vv{\lambda} = [\underline{\vv{\lambda}}^T,\bar{\vv{\lambda}}^T ] ^T \in\mathbb{R}^{2n}$ for optimization problem (\ref{eq:opt}), in which $\underline{\vv{\lambda}} = [ \underline{\lambda}_1,\ldots,\underline{\lambda}_n]^T$ and $\underline{\lambda}_i$ corresponds to the lower limit in the voltage constraint (\ref{eq:vol_con}); $\bar{\vv{\lambda}} = [ \bar{\lambda}_1,\ldots,\bar{\lambda}_n]^T$ and $\bar{\lambda}_i$ corresponds to the upper limit in (\ref{eq:vol_con}). We also introduce multiplier $\vv{\xi} = [\xi_1,\ldots,\xi_n]^T$ for the reactive power constraint (\ref{eq:q_con}), and $\xi_i$ is for both the lower limit and the upper limit in the capacity constraint (\ref{eq:q_con}). Next, we introduce the augmented Lagrangian for the optimization problem. {\color{black}Throughout this paper, we reserve letter $q$ for the physical reactive power injection, and use $\hatq$ to represent variables in the Lagrangian $\LL$.}
{\small
 \begin{align}
 \LL(\vv{\hatq},\vv{\xi},\vv{\lambda}) &= f(\vv{\hatq}) + \vv{\ulambda}^T(\vv{\underline{v}} - \vv{v}(\vv{\hatq})) + \vv{\blambda}^T( \vv{v}(\vv{\hatq}) - \vv{\bar{v}})+ K(\vv{\xi},\vv{\hatq}) \label{eq:lagrangian} 
 \end{align}}Here $K(\vv{\xi},\vv{\hatq}) = \sum_{i=1}^n K_i(\xi_i,\hatq_i)$, and $K_i(\xi_i,\hatq_i)$ is a quadratic penalty function defined to be, 
{\small
 \begin{align*}
 	K_i(\xi_i,\hatq_i) 
 	&= \left\{ \begin{array}{ll}
  \xi_i(\hatq_i - \unq_i) + \frac{c}{2}(\hatq_i - \unq_i)^2	 & \hatq_i +  \frac{\xi_i}{c}<\unq_i \\
  -\frac{\xi_i^2}{2c} &   \unq_i  \leq \hatq_i+   \frac{\xi_i}{c} \leq \barq_i  \\
  \xi_i(\hatq_i - \barq_i) + \frac{c}{2}(\hatq_i - \barq_i)^2	 & \hatq_i + \frac{\xi_i}{c} >\barq_i
   	 \end{array}  \right.
 \end{align*}}and we note that the partial derivatives of $K_i(\cdot,\cdot)$ are given as,
{\small\begin{align*}
 \frac{\partial K_i(\xi_i,\hatq_i)}{\partial \hatq_i}& = \text{ST}_{c \underline{q}_i}^{c\bar{q}_i}(\xi_i+ c \hatq_i)  \\
  \frac{\partial K_i(\xi_i,\hatq_i)}{\partial \xi_i} &= \frac{1}{c} [\text{ST}_{c \underline{q}_i}^{c\bar{q}_i}(\xi_i(t)+ c \hatq_i)   - \xi_i ].
\end{align*}}

In (\ref{eq:lagrangian}), $\vv{\ulambda}^T(\vv{\underline{v}} - \vv{v}(\vv{\hatq})) + \vv{\blambda}^T( \vv{v}(\vv{\hatq}) - \vv{\bar{v}})$ is the standard term in Lagrangian multiplier theory that penalizes violation of the voltage constraint, while the $K_i(\xi_i,\hatq_i)$ term is a special quadratic penalty function that penalizes violation of both the upper limit and the lower limit of constraint (\ref{eq:q_con}). For details of such quadratic penalty functions, we refer the readers to \cite[Section 3.2]{bertsekas2014constrained}, \cite[Appendix-G]{qu2018exponential}. In short, with the penalty functions, the max-min problem  
 \begin{equation}
 \max_{\vv{\lambda}\in\R^{2n}_{\geq 0}, \vv{\xi}\in\R^n} \min_{\vv{\hatq}\in\R^n} \LL(\vv{\hatq},\vv{\xi},\vv{\lambda}) \label{eq:maxmin}
 \end{equation}
 is equivalent to the original optimization problem (\ref{eq:opt}) (cf. Lemma~\ref{lem:saddle}). The reason we use the augmented Lagrangian instead of the standard Lagrangian is that the primal-dual gradient algorithm associated with the augmented Lagrangian avoids projection and has better convergence properties \cite{qu2018exponential}.  
 
We then write down the standard primal dual gradient algorithm \cite{nedic2009subgradient,cherukuri2016asymptotic,feijer2010stability,qu2018exponential} for solving the max-min problem (\ref{eq:maxmin}),
{\small \begin{subeqnarray} \label{eq:algo_concept}
 \hatq_i(t+1) & = &  \hatq_i(t) - \alpha  \frac{\partial  \LL(\vv{\hatq}(t), \vv{\xi}(t),\vv{\lambda}(t)) }{\partial \hatq_i} \nonumber \\
 &=& \hatq_i(t) - \alpha \Big[f_i'(\hatq_i(t)) + \text{ST}_{c \underline{q}_i}^{c\bar{q}_i}(\xi_i(t)+ c \hatq_i(t)) \nonumber  \\
 && + \sum_{j=1}^n [X]_{ij} (\bar{\lambda}_j(t) - \underline{\lambda}_j(t) + d \hatq_j(t))\Big]    \slabel{eq:algo_concept_1} \\
\xi_i(t+1) 
&= & \xi_i(t) + \beta \frac{\partial\LL(\vv{\hatq}(t), \vv{\xi}(t),\vv{\lambda}(t)}{\partial \xi_i } \nonumber \\
&=& \xi_i(t) + \beta \frac{\text{ST}_{c \underline{q}_i}^{c\bar{q}_i}(\xi_i(t)+ c \hatq_i(t))   - \xi_i}{c}  \slabel{eq:algo_concept_xi}
\\
 \blambda_i(t+1) & = & [ \blambda_i(t) + \gamma \frac{\partial \LL(\vv{\hatq}(t), \vv{\xi}(t),\vv{\lambda}(t) )}{\partial \blambda_i} ]^+ \nonumber\\
 &=& [\blambda_i(t) + \gamma (v_i(\vv{\hatq}(t)) - \bv_i) ]^+  \slabel{eq:algo_concept_blambda}\\
 \ulambda_i(t+1) & = & [\ulambda_i(t) + \gamma \frac{\partial \LL(\vv{\hatq}(t), \vv{\xi}(t),\vv{\lambda}(t) )}{\partial \ulambda_i}]^+ \nonumber\\
 &=& [\ulambda_i(t) + \gamma  (\uv_i - v_i(\vv{\hatq}(t))) ]^+  \slabel{eq:algo_concept_ulambda}
\end{subeqnarray} }

Eq. (\ref{eq:algo_concept}) is the standard primal-dual gradient algorithm (also known as saddle point algorithm). At every time step, $\vv{\hatq}(t)$ conducts a gradient \emph{descent} step along the gradient of $\LL$ w.r.t. $\vv{\hatq}$, since $\vv{\hatq}(t)$ seeks to minimize $\LL$ (cf. the max-min problem~\eqref{eq:maxmin}), while $\vv{\xi}(t)$ and $\vv{\lambda}(t)$ conducts a gradient \emph{ascent} step along the gradient of $\LL$ w.r.t. $\vv{\xi}$ and $\vv{\lambda}$, since $\vv{\xi}(t)$ and $\vv{\lambda}(t)$ seek to maximize $\LL$.
Though literature has shown the convergence of primal-dual gradient algorithms with properly chosen step sizes and under some conditions \cite{nedic2009subgradient,cherukuri2016asymptotic,feijer2010stability,qu2018exponential}, the algorithm in (\ref{eq:algo_concept}) does not meet our design requirements. Firstly, step (\ref{eq:algo_concept_1}) involves a summation from $1$ to $n$ and requires information across the network to implement, violating Requirement 1. Secondly, the $\hatq_i(t)$ in step (\ref{eq:algo_concept_1}) might violate the capacity constraint, violating Requirement 3. 

We now propose two modifications to (\ref{eq:algo_concept}) to meet the design requirements and the two modifications together change (\ref{eq:algo_concept}) into \algoname.

\textbf{Modification (a).} We now modify (\ref{eq:algo_concept_1}) such that each bus only needs local and neighbor's information to update. Eq. (\ref{eq:algo_concept_1}) is a gradient update for the $\vv{\hatq}$ coordinates of $\LL$, and the gradient is given by,
\begin{align}
&\nabla_{\vv{\hatq}} \LL (\vv{\hatq},\vv{\xi}(t),\vv{\lambda}(t))\nonumber\\
& = \nabla  f(\vv{\hatq}) +X(\vv{\bar{\lambda}}(t) - \vv{\underline{\lambda}}(t)) + \nabla_{\vv{\hatq}}K(\vv{\xi}(t),\vv{\hatq}).
\end{align}
Because of the sparse structure of $Y\coloneqq X^{-1}$ (cf. Proposition~\ref{prop:XY}), the scaled gradient $Y \nabla_{\vv{\hatq}}  \LL (\vv{\hatq},\vv{\xi}(t),\vv{\lambda}(t))$ is given by,
{\small\begin{align*}
[Y \nabla_{\vv{\hatq}}  \LL (\vv{\hatq},\vv{\xi}(t),\vv{\lambda}(t))]_i &=  \Big\{ \bar{\lambda}_i(t)  - \underline{\lambda}_i(t) +d \hatq_i \\
&\quad  + \sum_{j\in \mathcal{N}_i} [Y]_{ij} [f_j'(\hatq_j  ) + \text{ST}_{c\unq_j}^{c\barq_j}(\xi_j(t) + c \hatq_j )]  \Big\}
\end{align*}}
To calculate the $i$'th element of the scaled gradient $Y \nabla_{\vv{\hatq}}  \LL (\vv{\hatq},\vv{\xi}(t),\vv{\lambda}(t))$, bus $i$ only needs local information ($\blambda_i(t),\ulambda_i(t)$) and information from neighbors ($f'_j(\hatq_j), \text{ST}_{c\unq_j}^{c\barq_j}(\xi_j(t) + c \hatq_j )$ where $j\in\mathcal{N}_i$). Moreover, since $Y$ is positive definite, $Y \nabla_{\vv{\hatq}}  \LL (\vv{\hatq},\vv{\xi}(t),\vv{\lambda}(t))$ is still a descent direction for $\LL$ (in the $\vv{\hatq}$ coordinates) and hence using the scaled gradient in the primal dual gradient algorithm still has convergence guarantee \cite{cherukuri2016asymptotic,feijer2010stability}. Therefore, we change (\ref{eq:algo_concept_1}) into the following ``scaled'' gradient update, which gives rise to step (\ref{eq:algorithm_hatq}) in \algoname.
\begin{align}
\hatq_i(t+1) & =  \hatq_i(t) - \alpha [Y \nabla_{\vv{\hatq}}  \LL (\vv{\hatq}(t),\vv{\xi}(t),\vv{\lambda}(t))]_i \nonumber \\
&= \hat{q}_i(t) - \alpha  \Big\{ \bar{\lambda}_i(t)  - \underline{\lambda}_i(t) +d \hatq_i(t) \nonumber  \\
&\quad  + \sum_{j\in \mathcal{N}_i} [Y]_{ij} [f_j'( \hatq_j(t )  ) + \text{ST}_{c\unq_j}^{c\barq_j}(\xi_j(t) + c \hatq_j(t) )]  \Big\}. \label{eq:modify_distributed}
\end{align}


\textbf{Modification (b).} To fix the problem that $\hatq_i(t)$ may violate the capacity constraint, we do not actually implement $\hatq_i(t)$, but instead implement $q_i(t)  = [\hat{q}_i(t)]_{\unq_i}^{\barq_i}$, the projection of $\hatq_i(t)$ onto the capacity constraint. This gives rise to (\ref{eq:algorithm_q}) in our controller.  Another issue is that update (\ref{eq:algo_concept_blambda}) (\ref{eq:algo_concept_ulambda}) uses $v_i(\vv{\hatq}(t))$, which is not the measured voltage since the implemented reactive power is not $\vv{\hatq}(t)$. Therefore, we replace the $v_i(\vv{\hatq}(t))$ in (\ref{eq:algo_concept_blambda}) (\ref{eq:algo_concept_ulambda}) with $v_i(\vv{q}(t))$   and get,  
{\small\begin{align*}
 \blambda_i(t+1) & =  [ \blambda_i(t) + \gamma \frac{\partial \LL(\vv{q}(t),\vv{\xi}(t),\vv{\lambda}(t))}{\partial \blambda_i} ]^+ \nonumber\\
&= [\blambda_i(t) + \gamma (v_i(\vv{q}(t)) - \bv_i) ]^+  \slabel{eq:algo_concept_4}\\
\ulambda_i(t+1) & =  [\ulambda_i(t) + \gamma \frac{\partial \LL(\vv{q}(t),\vv{\xi}(t),\vv{\lambda}(t))}{\partial \ulambda_i}]^+ \nonumber\\
&= [\ulambda_i(t) + \gamma  (\uv_i - v_i(\vv{q}(t))) ]^+  \slabel{eq:algo_concept_5}
\end{align*}}which uses the measured voltage $v_i(\vv{q}(t)) = v_i(t)$ and gives rises to (\ref{eq:algorithm_blambda}) (\ref{eq:algorithm_ulambda}). We emphasize that after the change, the update for $\ulambda_i(t),\blambda_i(t)$ does not use the true gradient of $\LL$, $\nabla_{\vv{\lambda}}\LL(\vv{\hatq}(t),\vv{\xi}(t),\vv{\lambda}(t))  $ any more; but uses a gradient that is evaluated at a different point, $\nabla_{\vv{\lambda}}\LL(\vv{q}(t),\vv{\xi}(t),\vv{\lambda}(t)) $. We will show in Section \ref{sec:alg_analysis} that despite of the inconsistency, we will have $\vv{q}(t) - \vv{\hatq}(t) \rightarrow 0$ and the controller still converges. 

\subsection{Algorithm Analysis}\label{sec:alg_analysis}
We first analyze the Lagrangian (\ref{eq:lagrangian}), and show that the original problem (\ref{eq:opt}) is indeed equivalent to the max-min problem (\ref{eq:maxmin}). \ifthenelse{\boolean{isfullversion}}{The proof of Lemma~\ref{lem:saddle} is in Appendix-\ref{subsec:L}.}{The proof of Lemma~\ref{lem:saddle} can be found in \cite[Appendix-A]{fullversion}. }
\begin{lemma} \label{lem:saddle}
	$\LL(\vv{\hatq},\vv{\xi},\vv{\lambda})$  is convex in $\vv{\hatq}$, concave in $\vv{\xi}$, $\vv{\lambda}$ and
	 has a saddle point $(\vv{\hatq}_{sad},\vv{\xi}_{sad},\vv{\lambda}_{sad})$
	 satisfying $\LL(\vv{\hatq}_{sad},\vv{\xi}_{sad},\vv{\lambda}_{sad}) = \max_{\vv{\xi}\in\R^n,\vv{\lambda}\in\R^{2n}_{\geq 0}} \LL(\vv{\hatq}_{sad},\vv{\xi},\vv{\lambda}) = \min_{\vv{\hatq}\in\R^n} \LL(\vv{\hatq}, \vv{\xi}_{sad},\vv{\lambda}_{sad}) $. Moreover, for any saddle point $(\vv{\hatq}_{sad},\vv{\xi}_{sad},\vv{\lambda}_{sad})$, $\vv{\hatq}_{sad}$ must be the unique solution of the optimization problem~\eqref{eq:opt}. 
\end{lemma}

As discussed before, our algorithm is essentially the primal-dual gradient algorithm with scaled gradient, except that in the update for $\vv{\lambda}(t)$, the gradient is evaluated at $(\vv{q}(t),\vv{\xi}(t),\vv{\lambda}(t))$ instead of $(\vv{\hatq}(t), \vv{\xi}(t), \vv{\lambda}(t))$. Our proof essentially shows that $\vv{q}(t) - \vv{\hatq}(t) \rightarrow 0$, and therefore, our algorithm is approximately the true primal-dual gradient algorithm (with scaled gradient) and it converges to a saddle point of $\LL$. The rigorous proof of Theorem~\ref{thm:convergence} can be found \ifthenelse{\boolean{isfullversion}}{in Appendix-\ref{appendix:proof_convergence}.}{in our online report \cite[Appendix-B]{fullversion}.}
{\color{black}
\subsection{Implementation on a subset of nodes}\label{subsec:subcontrollable}
 Based on the derivations in Section~\ref{subsec:algo_derivation}, we can develop a variation of the algorithm \algoname\ that adapts to the case where only a subset of the buses are equipped with control components (referred to as controllable buses hereafter). Let the set of controllable buses be $\mathcal{C}$, and the set of uncontrollable buses be $\mathcal{U}$. Then, after rearranging the labels of the buses, we can write $X$ matrix as,
 $$ X = \left[\begin{array}{cc}
X_\mathcal{C} & G\\
G^T & X_\mathcal{U}
 \end{array}\right] $$  
and the voltage at the controllable buses $\vv{v}_\mathcal{C}$ is given by
 $$\vv{v}_\mathcal{C}(\vv{q}_\mathcal{C}) = X_\mathcal{C} \vv{q}_C + G \vv{q}_{\mathcal{U}} + \vv{v}_{\mathcal{U}}^{par}$$
We can formulate an analogous optimization problem as that in \eqref{eq:opt},\footnote{Similar to the $\frac{d}{2} q^T X q$ term in \eqref{eq:opt}, we can also add a term $\frac{d }{2} \vv{q}_\mathcal{C}^T X_\mathcal{C} \vv{q}_{\mathcal{C}}$ to \eqref{eq:opt_cost_sub} without changing too much of the discussion in this subsection. }
\begin{subeqnarray}\label{eq:opt_sub}
	\min_{\vv{q}_\mathcal{C}} & & f(\vv{q}_\mathcal{C}) \triangleq \sum_{i\in\mathcal{C}} f_i (q_i)  \slabel{eq:opt_cost_sub}\\
	s.t.  &  &\underline{v}_i \leq v_i(\vv{q}_{\mathcal{C}}) \leq \bar{v}_i, \forall i\in\mathcal{C} \slabel{eq:vol_con_sub}\\
	& &\underline{q}_i \leq q_i \leq \bar{q}_i, \forall i\in\mathcal{C} \slabel{eq:q_con_sub}
\end{subeqnarray}
where we note that the decision variables are now only $\vv{q}_\mathcal{C}$, the reactive power injection at the controllable buses, and the voltage constraint and the capacity constraint are also only placed on the controllable buses. Following the same modified primal-dual approach as in Section~\ref{subsec:algo_derivation}, we have the following resulting algorithm that operates on the controllable buses. For any $i\in\mathcal{C}$,

{\small \begin{subeqnarray} \label{eq:algo_sub}
		\hatq_i(t+1) 
		&=& \hatq_i(t) - \alpha \Big[ \bar{\lambda}_i(t)  - \underline{\lambda}_i(t)   \nonumber  \\
		&&  + \sum_{j\in \mathcal{C}} [X_\mathcal{C}^{-1}]_{ij} [f_j'( \hatq_j(t )  ) + \text{ST}_{c\unq_j}^{c\barq_j}(\xi_j(t) + c \hatq_j(t) )]  \Big]  \nonumber\\
		  \slabel{eq:algo_sub_1} \\
		\xi_i(t+1) 
		&=& \xi_i(t) + \beta \frac{\text{ST}_{c \underline{q}_i}^{c\bar{q}_i}(\xi_i(t)+ c \hatq_i(t))   - \xi_i}{c}  \slabel{eq:algo_sub_xi}
		\\
		\blambda_i(t+1) 
		&=& [\blambda_i(t) + \gamma (v_i(t) - \bv_i) ]^+  \slabel{eq:algo_sub_blambda}\\
		\ulambda_i(t+1) 
		&=& [\ulambda_i(t) + \gamma  (\uv_i - v_i(t)) ]^+  \slabel{eq:algo_sub_ulambda}\\
		  q_i(t+1) &=& [\hat{q}_i(t+1)]_{\unq_i}^{\barq_i}  \slabel{eq:algo_sub_q}
\end{subeqnarray} }The same theoretic guarantee holds for the above algorithm \eqref{eq:algo_sub}. The only substantial difference is that the matrix $X_\mathcal{C}^{-1}$ has a sparsity pattern different from that of the physical network topology. In detail, $i$ and $j$ can communicate if and only if the path between $i$ and $j$ in the physical network does not pass any other bus in $\mathcal{C}$. An illustration of this communication topology is given in Figure~\ref{fig:comm_novar}. For more information, we also refer to \cite[Definition 4]{cavraro2016value}. \textcolor{black}{This local communication still has advantage over a centralized scheme where all the buses need to communicate with a central point in presence of communication challenges such as delay, limited bandwidth, node failures, etc.}
\begin{figure}
	\centering
	\includegraphics[width=0.7\columnwidth]{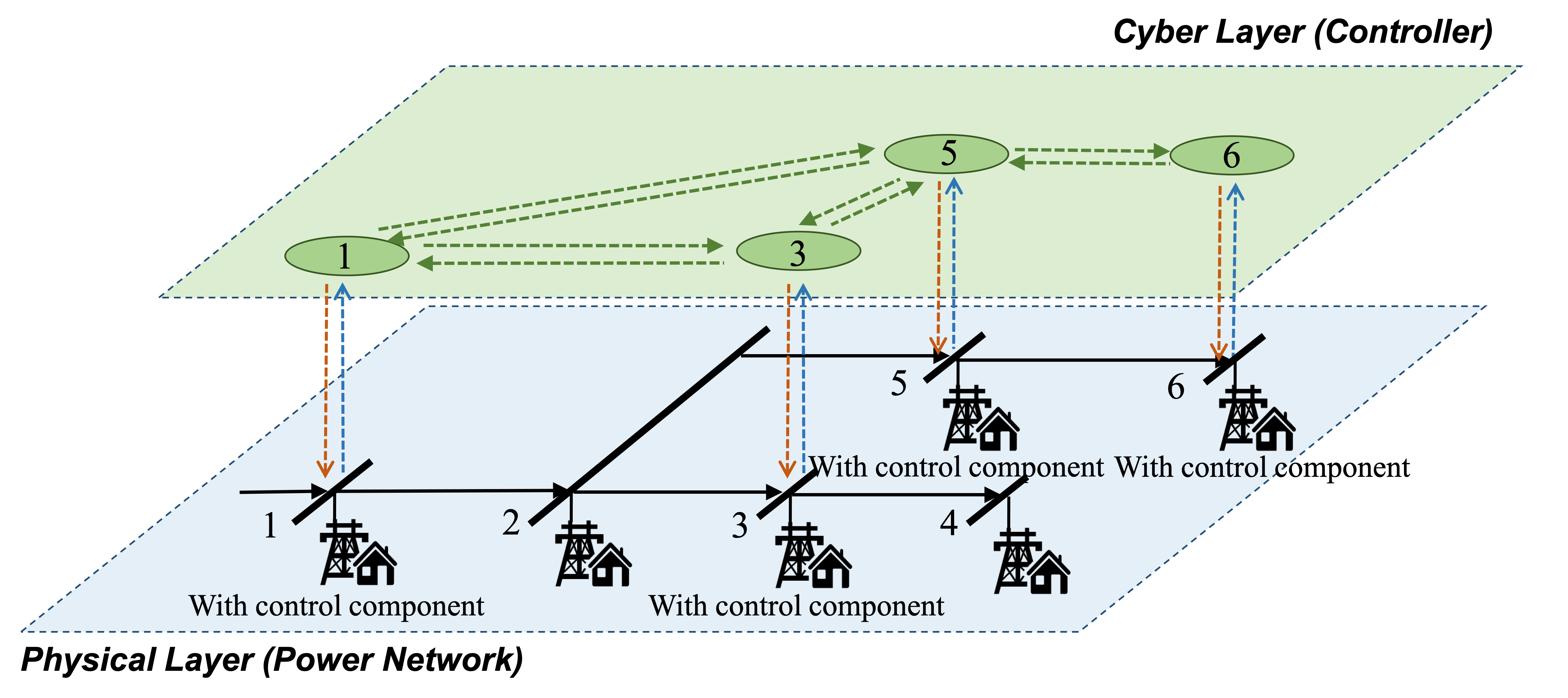}
	\caption{Communication graph when some buses do not have control components. In the above network, buses $\mathcal{C}=\{1,3,5,6\}$ are equipped with control components while buses $\mathcal{U} = \{2,4\}$ do not. In the communication graph, bus $1$ can communicate with bus $3$ because the path between $1$ and $3$ passes through $2$, a bus \emph{not} in $\mathcal{C}$. Bus $1$ cannot communicate with bus $6$ as the path between $1$ and $6$ passes through bus $2,5$, and bus $5$ is in $\mathcal{C}$. }\label{fig:comm_novar}
\end{figure}
}

\section{case study}\label{sec:case}

\subsection{Single Phase Network}
 \begin{figure}[!t]
 	\centering
 	\includegraphics[width=0.7\columnwidth]{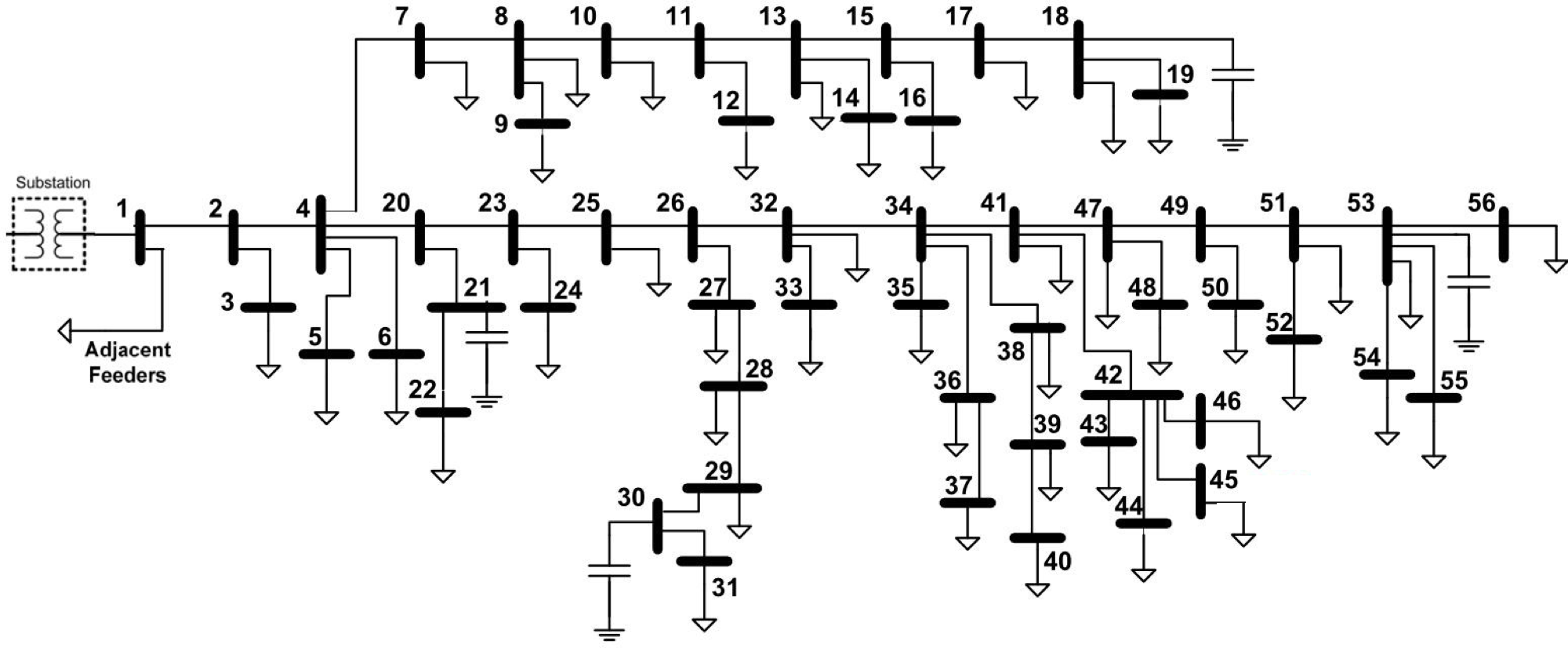}
 	\caption{Schematic diagram of two SCE distribution systems. }
 	\label{fig:circuit}
 \end{figure}
We evaluate \algoname\ on a distribution circuit of South California Edison with a high penetration of photovoltaic (PV) generation \cite{Farivar-2012-VVC-PES}.
Figure \ref{fig:circuit} shows the $56$-bus distribution circuit.  
Note that bus $1$ indicates the substation, and there is PV generation at various locations of the network (bus $2, 4, 7,8, 9,10,11,12,13,14,15,16,19,20,23,25,26,32$). See \cite{Farivar-2012-VVC-PES} for the network data. 

In the simulation, we assume that there are control components at all the buses and those control components can supply or consume at most $0.2$ MVar reactive power (i.e. $\barq_i = 0.2 MVar, \unq_i = -0.2MVar$). The nominal voltage magnitude is $12 \text{kV}$ and the acceptable range is set as $[11.4\text{kV}, 12.6\text{kV}]$ which is the plus/minus 5\% of the nominal value. 
{\color{black}In line with Remark~\ref{rem:cost}, we set $f_i (q_i) =  \frac{\eta}{s_i^{max}} q_i^2$ where $\eta=10$ and $s_i^{max}$ are synthetic values lying between $0.5$ and $1$. We set the parameter $d = 1 $.}
Though the analysis of this paper is built on the linearized power flow model (\ref{eq:branchflow_linear}),  we simulate \algoname\ with the full nonlinear AC power flow model (\ref{eq:branchflow_AC}) using MATPOWER \cite{5491276}.

\begin{figure*}[t]
	\begin{subfigure}[b]{0.33\textwidth}
		\centering
		\includegraphics[width=\textwidth]{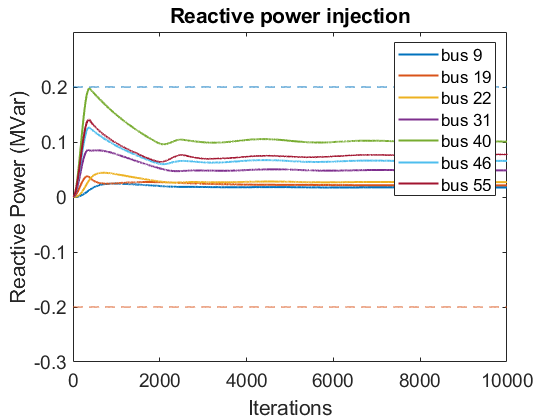} 
	\end{subfigure}
	\begin{subfigure}[b]{0.33\textwidth}
	\centering
	\includegraphics[width=\textwidth]{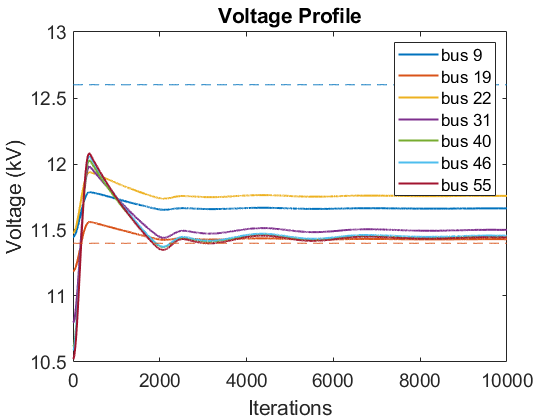} 
\end{subfigure}
	\begin{subfigure}[b]{0.33\textwidth}
		\centering
		\includegraphics[width=\textwidth]{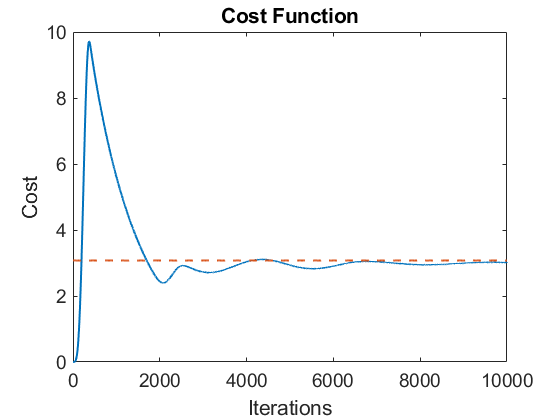} 
	\end{subfigure}
	\caption{Simulation results for the case with static load and PV generation. The left plot and the middle plot are the reactive power injection (control action) and the voltage profile for a selected subset of buses. The right plot is the cost function.}\label{fig:static}
\end{figure*}

{\color{black}
\textit{Static load and PV generation.} Firstly, we run \algoname\ in a scenario that the load and the PV generation are not time-varying. We consider a heavily loaded setting, resulting in low voltages at the buses. The simulation results are given in Figure \ref{fig:static}. It shows that \algoname\ can bring the voltage to the acceptable range and in the meanwhile not violating the capacity constraint. Further, the dashed line in the right plot of Figure \ref{fig:static} is the optimal value of the optimization problem \eqref{eq:opt} with the AC power flow model \eqref{eq:branchflow_AC}, which we calculate through the SOCP (Second Order Cone Programming) convex-relaxation method in \cite{low2014convex}. The right plot of Figure \ref{fig:static} shows that \algoname\ is able to drive the system to a near-optimal operating point.  }

\begin{figure}
	\begin{center}
	\includegraphics[width = 0.5\columnwidth]{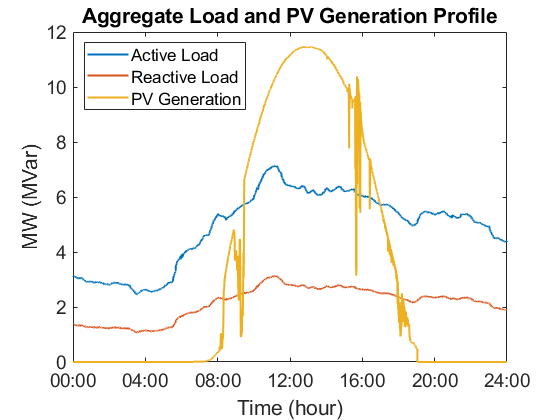}
		\end{center}
	\caption{Aggregate load and PV generation profile.}\label{fig:load}
\end{figure}

\begin{figure*}[t]

	\begin{subfigure}[b]{0.33\textwidth}
		\centering
		\includegraphics[width=\textwidth]{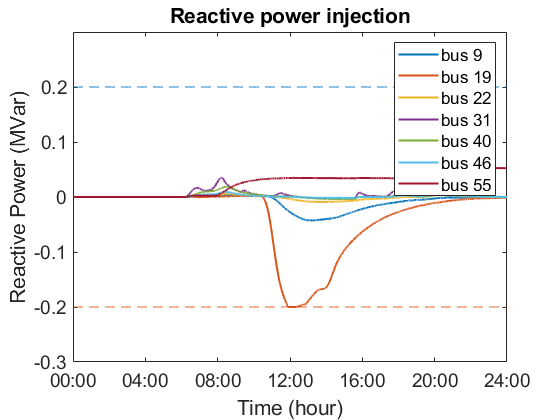} 
	\end{subfigure}
	\begin{subfigure}[b]{0.33\textwidth}
	\centering
	\includegraphics[width=\textwidth]{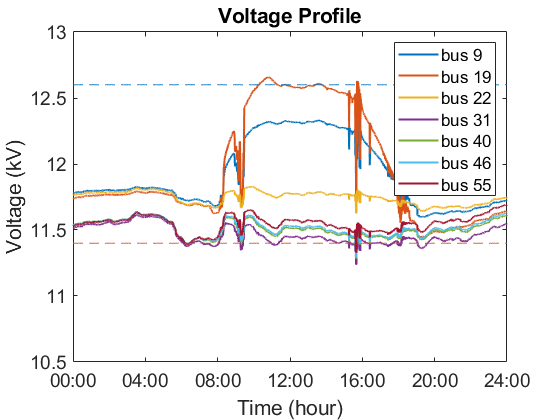} 
\end{subfigure}
	\begin{subfigure}[b]{0.33\textwidth}
		\centering
		\includegraphics[width=\textwidth]{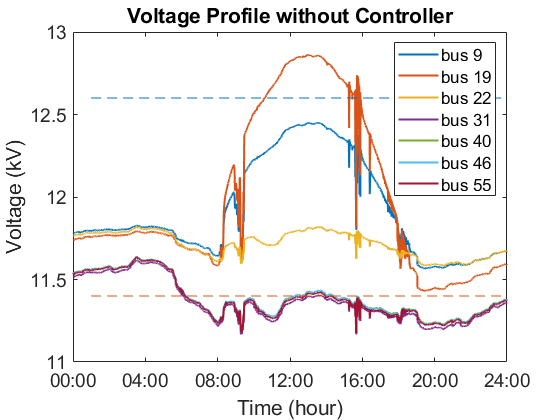} 
	\end{subfigure}
	\caption{Simulation results for the case with time-varying load and PV generation. The left and the middle plot are the reactive power injection (control action) and the voltage profile for a selected subset of buses when implementing \algoname. The right plot is the voltage profile when no voltage controller is used. }\label{fig:time_varying}
\end{figure*}

\textit{Time-varying load and PV generation.} Next, we test \algoname\ in a more realistic setting. {We use the load and PV generation profile  in \cite{bernstein2017real}.} The time span of the data set is one day (24 hours), and the time resolution is $6$s. We plot the aggregated load and PV generation profile across the buses in Figure~\ref{fig:load}. Figure~\ref{fig:load} shows that the load increases after approximately 6AM; the PV generation is.nonzero between 8AM and 7PM, peaks at noon, and has large fluctuations throughout the day. Consistent with the time resolution of the dataset, we identify each iteration in \algoname\ with $6$ seconds, which means that the controller adjusts its control action every $6$ seconds. We run \algoname\ in this setting and simulate the voltage profile and the reactive power injection. For comparison, we also simulate the network voltage profile when no voltage controller is used. The simulation results are given in Figure \ref{fig:time_varying}. It shows that, despite the volatility in load and PV generation, \algoname\ can quickly bring the voltage into the acceptable range and in the mean while, not violating the capacity constraint.

\begin{figure*}[t]
	\begin{subfigure}[t]{0.33\textwidth}
		\centering
		\includegraphics[width=\textwidth]{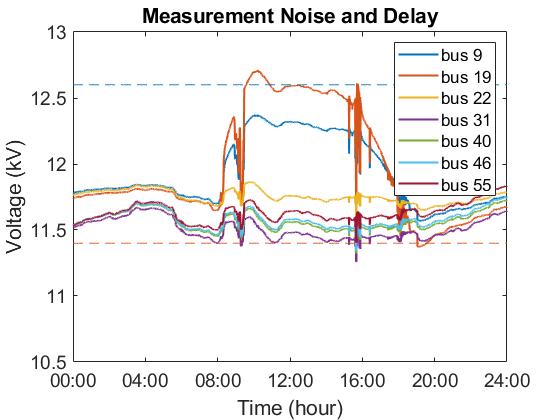} 
		\caption{Robustness against measurement noise and delay.}\label{fig:robust_noise}
	\end{subfigure}
	\begin{subfigure}[t]{0.33\textwidth}
		\centering
		\includegraphics[width=\textwidth]{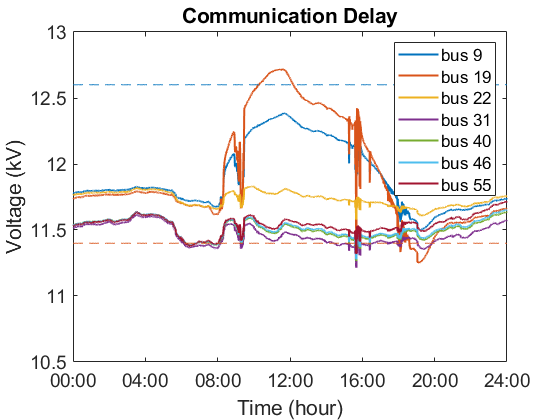} 
				\caption{Robustness against communication delay.}\label{fig:robust_delay}
	\end{subfigure}
	\begin{subfigure}[t]{0.33\textwidth}
		\centering
		\includegraphics[width=\textwidth]{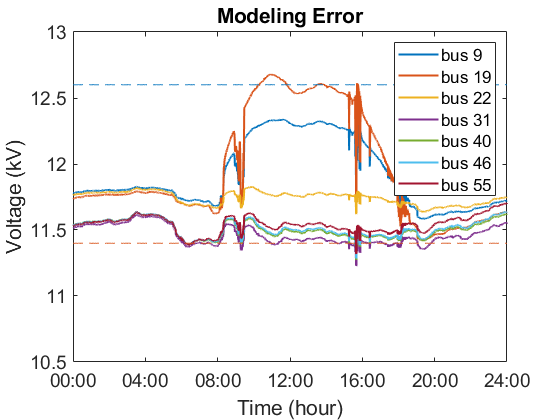} 
						\caption{Robustness against modeling error.}\label{fig:robust_model}
	\end{subfigure}
	\caption{Robustness of \algoname\ against measurement noise and delay, communication delay and modeling error.}\label{fig:robust}
\end{figure*}

{\color{black}\textit{Impact of measurement error and delay.} We test the robustness of \algoname\ against measurement error and delay. We use the same simulation setting as the time-varying load and PV generation case (Figure~\ref{fig:time_varying}), except that each voltage measurement is corrupted by a random Gaussian noise with zero mean and $0.03$ p.u. ($0.36$kV) standard deviation, and delayed by $30s$ ($5$ time steps). The result is shown in Figure~\ref{fig:robust_noise}. It can be seen that under noisy and delayed measurements, \algoname\ can still guarantee that the voltage lies within the limits for most of the time, though compared to Figure~\ref{fig:time_varying} there are larger voltage violations between 8:00 to 12:00.}

{\color{black}\textit{Impact of communication delay.} We test the robustness of \algoname\ against communication delay. In particular, when node $j$ sends variables to node $i$, we let the communication be delayed for a random amount of time up to $10$ iterations ($60$ seconds). This means that the controller at bus $i$ implements step \eqref{eq:algorithm_hatq} in the following way
	
{\small	\begin{align*}
	 	\hat{q}_i(t+1)  &=  \hat{q}_i(t) - \alpha  \bigg\{ \bar{\lambda}_i(t)  - \underline{\lambda}_i(t) + d \hatq_i(t)  \nonumber \\
	&\quad  + \sum_{j\in \mathcal{ N}_i} [Y]_{ij} \Big[f_j'(\hat{q}_j(t - \Delta_{ij}^t )  ) \\
	&\quad+ \text{ST}_{c \underline{q}_j}^{c\bar{q}_j}(\xi_j(t- \Delta_{ij}^t)+ c \hatq_j(t- \Delta_{ij}^t)) \Big]  \bigg\}
	\end{align*}}where the information received from $j$ is delayed by $\Delta_{ij}^t$ steps, and $\Delta_{ij}^t$ is a random integer between $0$ and $10$. We implement \algoname\ under the above setting, and plot the results in Figure~\ref{fig:robust_delay}. It can be seen that in the delayed communication case \algoname\ can still guarantee the voltage stays inside the acceptable range for most of the time. Compared with Figure~\ref{fig:time_varying}, the voltage trajectory under communication delay shown in Figure~\ref{fig:robust_delay}  violates the voltage constraint by larger amounts. }

\textit{Impact of modeling error.} Recall that the implementation of \algoname\ requires knowledge of $[Y]_{ij}$, which is essentially the reactance of the power lines (cf. Proposition~\ref{prop:XY}). We test \algoname\ with incorrect values of $[Y]_{ij}$ that are within plus/minus $20\%$ of the true value. The results are shown in Figure~\ref{fig:robust_model}. Figure~\ref{fig:robust_model} shows that under inaccurate model \algoname\ performs similarly as the case with accurate model (Figure~\ref{fig:time_varying}).

{\color{black}\subsection{Three-Phase Unbalanced Network}
We evaluate \algoname\ on a three-phase unbalanced network.\footnote{The three phase simulation is based on the models in \cite[Section V]{bernstein2018load}. Since high fidelity three phase models that are suited for optimization and control algorithm design are still being actively developed, our future work includes rerun these tests using newly developed models and softwares.} We use the network with 126 multiphase nodes, the PV generation and the load profile in \cite{bernstein2017real}. See \cite[Section V]{bernstein2017real} for a detailed description of the network structure and parameters. The control components are located at parts of the buses of the network, and for the buses that have control components, the control components are located at parts of the phases. We run \algoname\ on the network, and show the voltage profile of phase a of the buses that have control components in Figure~\ref{fig:3phase} (left). The voltage profile of the other phases is similar. We also show the voltage profile when no controller is used in Figure~\ref{fig:3phase} (right). Figure~\ref{fig:3phase} shows that \algoname\ can maintain the voltage inside the acceptable limit despite the network is unbalanced.
\begin{figure*}[t]
	\begin{subfigure}[t]{0.5\textwidth}
		\centering
		\includegraphics[width=\textwidth]{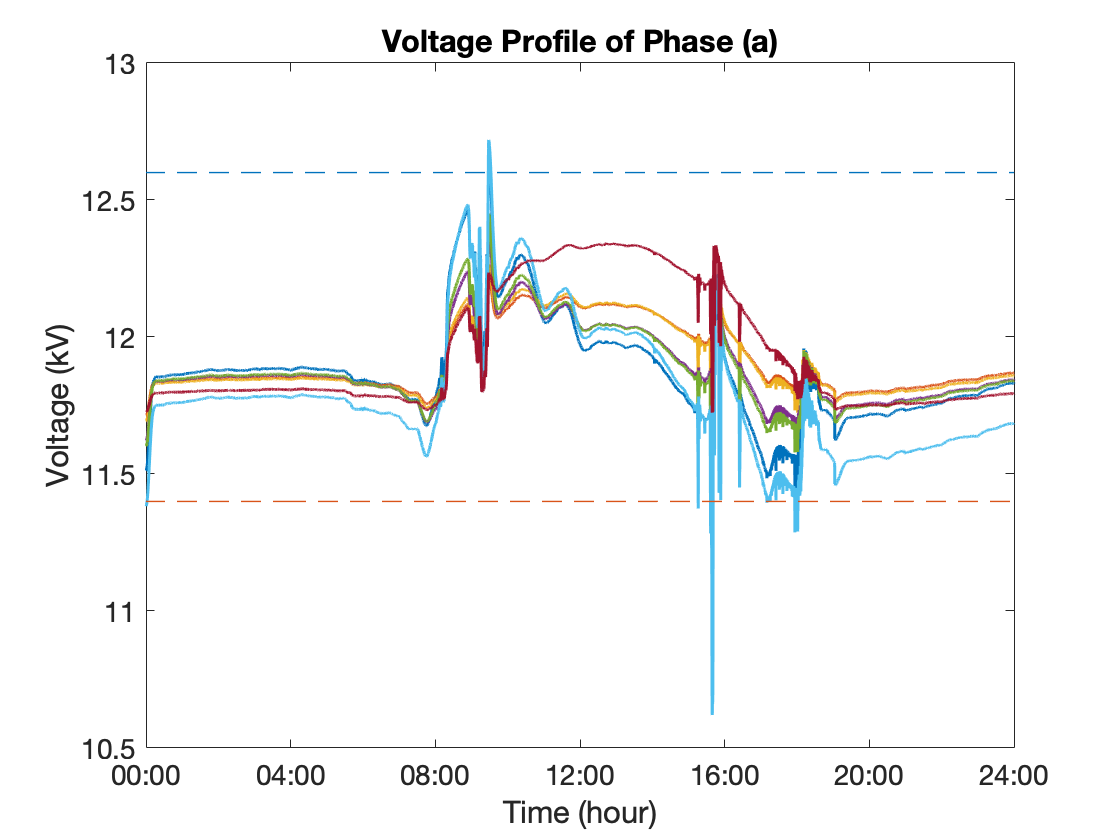} 
	\end{subfigure}
	\begin{subfigure}[t]{0.5\textwidth}
		\centering
		\includegraphics[width=\textwidth]{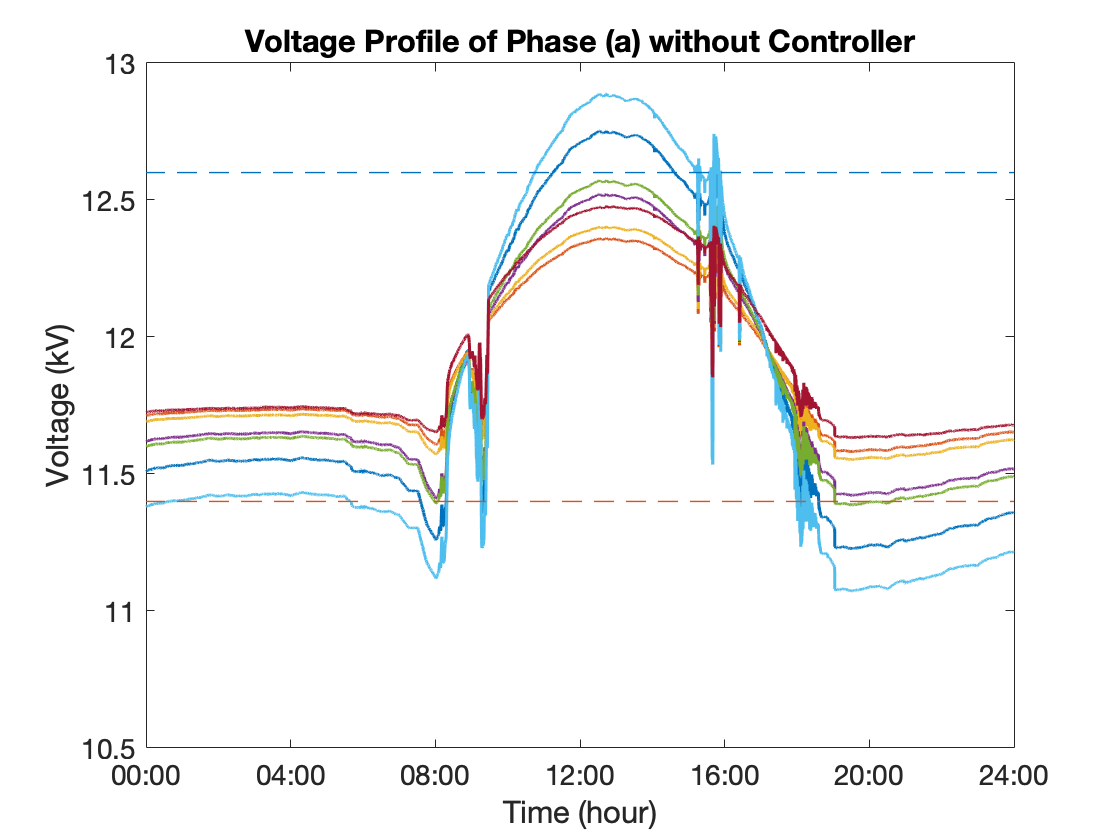} 
	\end{subfigure}
	\caption{Simulation results for three-phase unbalanced network. Left: when \algoname\ is used. Right: when no control is used. }\label{fig:3phase}
\end{figure*}}

\section{Conclusion}\label{sec:conclusion}
This paper proposes a distributed feedback voltage controller \algoname\ that can (i) meet the voltage constraint asymptotically, (ii) satisfy the reactive power capacity constraint throughout, and (iii) minimize a cost that is composed of a power loss related cost and reactive power operation costs. Future work includes extending the approach to jointly control active and reactive power.

\bibliographystyle{IEEETran}
\bibliography{VoltageRef-2017}
\ifthenelse{\boolean{isfullversion}}{

\appendix 
\subsection{Proof of Lemma~\ref{lem:saddle}}\label{subsec:L}


For notational simplicity, we define $\tilde{X} = [-X, X]^T$, $\vv{v}_b = [-\vv{ \uv}^T + (\vv{v}^{par})^T, \vv{\bv}^T -(\vv{v}^{par})^T ]^T$ and rewrite Lagrangian (\ref{eq:lagrangian}) as,
\begin{align}
\LL(\vv{\hatq},\vv{\xi},\vv{\lambda})  & = f(\vv{\hatq}) + \vv{\lambda}^T ( \tilde{X} \vv{\hatq} - \vv{v}_b)    + K(\vv{\xi},\vv{\hatq}).
\end{align}

	 Recall that $(\vv{\hatq},\vv{\xi},\vv{\lambda})$ is a saddle point if and only if $\LL(\vv{\hatq},\vv{\xi},\vv{\lambda}) = \min_{\vv{\hatq}'\in\R^n} \LL(\vv{\hatq}',\vv{\xi},\vv{\lambda}) =\max_{\vv{\xi}'\in\R^n}\LL(\vv{\hatq},\vv{\xi}',\vv{\lambda})= \max_{\vv{\lambda}'\in\R^{2n}_{\geq 0}}\LL(\vv{\hatq},\vv{\xi},\vv{\lambda}')$. Notice that $\LL(\vv{\hatq},\vv{\xi},\vv{\lambda})$ is convex and lower bounded in $\vv{\hatq}$, concave and upper bounded in $\vv{\xi}$, and linear in $\vv{\lambda}$, so $(\vv{\hatq},\vv{\xi},\vv{\lambda})$ is a saddle point if and only if $\forall i$,
	 \begin{align}
	 \frac{\partial}{\partial \hatq_i} \LL(\vv{\hatq},\vv{\xi},\vv{\lambda}) &= 0 \\
	 \frac{\partial}{\partial\xi_i} \LL(\vv{\hatq},\vv{\xi},\vv{\lambda}) & = 0 \\
\ulambda_i	 \frac{\partial}{\partial\ulambda_i} \LL(\vv{\hatq},\vv{\xi},\vv{\lambda}) =	\bar{\lambda}_i \frac{\partial}{\partial\blambda_i} \LL(\vv{\hatq},\vv{\xi},\vv{\lambda})  &=0 \label{eq:saddle:lambda}\\
 \ulambda_i \geq 0,	\blambda_i &\geq 0  \label{eq:saddle:dual_feasibility_lambda} \\
\uv_i\leq v_i(\vv{\hatq})&\leq\bv_i. \label{eq:saddle:primal_feasibility}
	 \end{align}
We claim the above equations are equivalent to the KKT condition of optimization problem \eqref{eq:opt}. First notice that $ \frac{\partial}{\partial\xi_i} \LL(\vv{\hatq},\vv{\xi},\vv{\lambda})  = 0$ is 
	$$   \text{ST}_{c\unq_i}^{c\barq_i} (\xi_i + c\hatq_i) = \xi_i,\forall i. $$
	Recall that $\text{ST}_{c\unq_i}^{c\barq_i}(\cdot)$ is the soft thresholding function, and we can check that the above equation is equivalent to
	\begin{align}
	  \hatq_i & \geq \unq_i \label{eq:saddle:pf_1}\\
	  \hatq_i & \leq \barq_i \label{eq:saddle:pf_2}\\
	  (\hatq_i - \barq_i)\max(\xi_i,0)&=0 \label{eq:saddle:cs_1}\\
	  (\hatq_i - \unq_i)\min(\xi_i,0)&= 0  \label{eq:saddle:cs_2}
	\end{align}
	Next, notice that \eqref{eq:saddle:lambda} is equivalent to, 
	\begin{align}
	\ulambda_i(v_i(\vv{\hatq}) - \uv_i) = \blambda_i(v_i(\vv{\hatq}) - \bv_i) = 0.\label{eq:saddle:cs_3}
	\end{align}
	Further, $\nabla_{\hatq} \LL(\hatq,\xi,\lambda) = 0 $ can be rewritten as,
{\small	\begin{align}
	0 &= \nabla f(\vv{\hatq}) + \tilde{X}^T \vv{\lambda} + \sum_{i=1}^n \text{ST}_{c\unq_i}^{c\barq_i} (\xi_i + c\hatq_i) \vv{e}_i \nonumber \\
	&= \nabla f(\vv{\hatq}) + \tilde{X}^T \vv{\lambda} + \sum_{i=1}^n \xi_i \vv{e}_i\nonumber \\
	&= \nabla f(\vv{\hatq}) + \tilde{X}^T \vv{\lambda} + \sum_{i=1}^n \max(\xi_i,0) \vv{e}_i + \sum_{i=1}^n (-\min(\xi_i,0)) (-\vv{e}_i)  \label{eq:saddleinner:stationarity}
	\end{align}}where $\vv{e}_i\in\R^n$ is a vector where the $i$'th entry is $1$ and other entries are $0$.

Now we replace $\max(\xi_i,0) $ with $\bar{\xi}_i$, and replace $-\min(\xi_i,0)$ with $\underline{\xi}_i$ and impose constraint $\bar{\xi}_i\geq 0$, $\underline{\xi}_i\geq 0$.\footnote{It is easy to check such replacement is one-to-one and onto. } Then \eqref{eq:saddle:dual_feasibility_lambda}\eqref{eq:saddle:primal_feasibility}\eqref{eq:saddle:pf_1}\eqref{eq:saddle:pf_2}\eqref{eq:saddle:cs_1}\eqref{eq:saddle:cs_2}\eqref{eq:saddle:cs_3}\eqref{eq:saddleinner:stationarity}, together with $\bar{\xi}_i\geq 0$, $\underline{\xi}_i\geq 0$, are exactly the KKT condition of optimization problem~\eqref{eq:opt}.\footnote{\eqref{eq:saddleinner:stationarity} is the stationarity condition, \eqref{eq:saddle:primal_feasibility}\eqref{eq:saddle:pf_1}\eqref{eq:saddle:pf_2} are the primal feasibility condition, \eqref{eq:saddle:dual_feasibility_lambda} and $\bar{\xi}_i\geq 0,\underline{\xi}_i\geq 0$ are the dual feasibility condition, \eqref{eq:saddle:cs_1}\eqref{eq:saddle:cs_2}\eqref{eq:saddle:cs_3} are the complimentary slackness condition. } 

Since the slater condition of problem \eqref{eq:opt} holds, we have the KKT condition of problem \eqref{eq:opt} has a solution, and hence $\LL$ has a saddle point $(\vv{\hatq}_{sad},\vv{\xi}_{sad},\vv{\lambda}_{sad})$. Moreover, since the objective of \eqref{eq:opt} is strongly convex, the primal solution is unique. Therefore, for any saddle point $(\vv{\hatq}_{sad},\vv{\xi}_{sad},\vv{\lambda}_{sad})$ of $\LL(\cdot,\cdot,\cdot)$, $\vv{\hatq}_{sad}$ must be the unique solution of problem \eqref{eq:opt}. 
\qedd

\subsection{Proof of Theorem~\ref{thm:convergence}}\label{appendix:proof_convergence}
As mentioned in Section \ref{subsec:concept}, \algoname\ is not the standard primal dual gradient algorithm because in update (\ref{eq:algorithm_blambda}) (\ref{eq:algorithm_ulambda}), the gradient of $\LL$ is not evaluated at $\vv{\hatq}(t)$, but at $\vv{q}(t)$ instead. In what follows, we define
\begin{equation}
\mathcal{S}(\vv{\lambda}) = \max_{\vv{\xi}}\min_{\vv{\hatq}} \LL(\vv{\hatq},\vv{\xi},\vv{\lambda})
\end{equation}
i.e. $ \mathcal{S}(\vv{\lambda}) $ is the solution of the inner layer of the max-min problem~\eqref{eq:maxmin} given fixed $\vv{\lambda}$. 
Our analysis treats the update of $\vv{\lambda}(t)$ as the gradient projection algorithm for $S(\vv{\lambda})$ where the gradient contains some ``error''. To fully understand the property of $\SSS$, we prove the following Lemma on the max-min problem of $\LL(\cdot,\cdot,\vv{\lambda})$ given fixed $\vv{\lambda}$. The proof of Lemma~\ref{lem:L_property} is given in Appendix-\ref{appendix:proof_L}.


\begin{lemma} \label{lem:L_property} The following statements hold.
	\begin{itemize} 
		\item [(a)] For every $\vv{\lambda}$, function $\LL(\cdot,\cdot,\vv{\lambda})$ has a unique saddle point $(\vv{\hatq}^*(\vv{\lambda}),\vv{\xi}^*(\vv{\lambda}))$ satisfying $\LL(\vv{\hatq}^*(\vv{\lambda}),\vv{\xi}^*(\vv{\lambda}),\vv{\lambda}) =\min_{\vv{\hatq}\in\R^n}\LL(\vv{\hatq},\vv{\xi}^*(\vv{\lambda}),\vv{\lambda}) = \max_{\vv{\xi}\in\R^n}\LL(\vv{\hatq}^*(\vv{\lambda}),\vv{\xi},\vv{\lambda})$. Moreover, $\hatq^*_i(\vv{\lambda}) \in[\unq_i,\barq_i]$. 
		\item [(b)] $\vv{\xi}^*(\vv{\lambda})$ and $\vv{\hatq}^*(\vv{\lambda})$ is Lipschitz in $\vv{\lambda}$.
		\begin{align}
		\Vert \vv{\hatq}^*(\vv{\lambda}) - \vv{\hatq}^*(\vv{\lambda}')\Vert & \leq \frac{\Vert\tildex\Vert}{\mu}\Vert\vv{\lambda}-\vv{\lambda}'\Vert  \\
		\Vert \vv{\xi}^*(\vv{\lambda}) - \vv{\xi}^*(\vv{\lambda}') \Vert &\leq 2 \frac{l}{\mu}\Vert \tildex\Vert\Vert\vv{\lambda}-\vv{\lambda}'\Vert
		\end{align}
		\item [(c)] $\SSS(\vv{\lambda})$ is a concave function, and $\nabla \SSS(\vv{\lambda}) = \tildex \vv{\hatq}^*(\vv{\lambda}) -\vv{v}_b$, and $\SSS$ is $\frac{\Vert\tildex\Vert^2}{\mu } $-smooth. 
		\item [(d)] $\SSS(\vv{\lambda})$ is upper bounded over $\vv{\lambda}\in \R^{2n}_{\geq 0}$ and moreover, for any real number $a\in\R$, level set $Level_a = \{\vv{\lambda} \in\R^{2n}_{\geq 0}:\SSS(\vv{\lambda}) \geq a \}$ is bounded. 
		\item [(e)] Let $\vv{\lambda}^*$ be any solution of $\max_{\vv{\lambda}\in\R^{2n}_{\geq 0}} \SSS(\vv{\lambda})$, then $\vv{\hatq}^*(\vv{\lambda}^*)$ is the \textit{unique} solution of optimization problem (\ref{eq:opt}).  
	\end{itemize}
\end{lemma}
Now we show that the update of $\vv{\lambda}(t)$ can be written as the gradient projection algorithm for $S(\vv{\lambda})$ with noisy gradient. For notational simplicity, we abuse the notation $\vv{\hatq}^*(\cdot), \vv{\xi}^*(\cdot)$ in Lemma~\ref{lem:L_property} and define $\vv{\hatq}^*(t) := \vv{\hatq}^*(\vv{\lambda}(t)), \vv{\xi}^*(t) := \vv{\xi}^*(\vv{\lambda}(t)).$ Then, the gradient of $\mathcal{S}$ at $\vv{\lambda}(t)$ is (by Lemma~\ref{lem:L_property} (c)) 
$$\nabla \mathcal{S} (\vv{\lambda}(t))  = \tildex \vv{\hatq}^*(t) - \vv{v}_b.$$
In the meanwhile, notice the update for $\vv{\lambda}(t)$ (\ref{eq:algorithm_blambda}) (\ref{eq:algorithm_ulambda}) can be written as, 
$$\vv{\lambda}(t+1)=[\vv{\lambda}(t) + \gamma (\tildex \vv{q}(t) - \vv{v}_b)]^+ .$$
So the update for $\vv{\lambda}(t)$ is ``inexact'' projected gradient method for function $\SSS$,  and the inexact gradient is $ \tildex \vv{q}(t) - \vv{v}_b$. Therefore, the gradient error is $$\epsilon(t) = \Vert \tildex (\vv{q}(t) - \vv{\hatq}^{*}(t))\Vert. $$

We show in Lemma~\ref{lem:grad_error} the gradient error is bounded, whose proof is deferred to Appendix-\ref{subsec:bound_gd_error}. 
\begin{lemma}\label{lem:grad_error}
	The gradient error $\epsilon(t) = \Vert  \tildex (\vv{q}(t) - \vv{\hatq}^{*}(t)) \Vert$ is bounded by
	\begin{align*}
	\epsilon(t)&\leq \Vert \tildex\Vert \Vert \vv{q}(t) - \vv{\hatq}^*(t)\Vert \leq \Vert \tildex\Vert \Vert \vv{\hatq}(t) - \vv{\hatq}^*(t)\Vert \\
	& \leq C_1 \rho^t + C_2 \sum_{k=0}^{t-1} \rho^{t-1-k} \Vert \vv{\lambda}(k+1) - \vv{\lambda}(k)\Vert 
	\end{align*}
	for constants $\rho$, $C_1$ and $C_2$ defined as follows. First, define $\mu'=\mu \underline{\sigma}(Y),l' = l\bar{\sigma}(Y)$. Further define constants,
	\begin{align*}
	a &= 20 l' [\max( \frac{c \bar{\sigma}(Y)}{\mu'}, \frac{l'}{\mu'} )]^2 [\max(\frac{\beta}{\alpha l' c},\frac{l'}{\mu'} ) ]^2 \kappa(Y)\\
	\tau& = \frac{\beta }{2\alpha } \frac{\underline{\sigma}(Y)}{a}\\
	\nu &= \max( \sqrt{2(l'+ \bar{\sigma}(Y)c)^2 + 2 \frac{\beta^2}{\alpha^2} \bar{\sigma}(Y)},\sqrt{2 \bar{\sigma}(Y) + 8 \frac{\beta^2}{\alpha^2 c^2}  })\\
P_0 &= \left[\begin{array}{cc}
	\frac{\beta}{\alpha}a  I         & \frac{\beta}{\alpha} Y^{1/2}\\
	\frac{\beta}{\alpha} Y^{1/2} & a I
	\end{array}\right]\\
	 P& = \left[\begin{array}{cc}
	\frac{\beta}{\alpha} a Y^{-1}         & \frac{\beta}{\alpha} I\\
	\frac{\beta}{\alpha} I  & a I
	\end{array}\right]
			\end{align*}
			Then, we can define $\rho$, $C_1$, and $C_2$ as follows,
			\begin{align*}
 \rho &= (e^{-\frac{\tau\alpha}{2}} + \frac{\alpha^2 \nu^2 \kappa(P_0) }{2} )	\\		
 C_1 &= \Vert\tildex \Vert\sqrt{\kappa(P)} \Vert \vv{z}(0) - \vv{z}^*(0)\Vert \\
 C_2 &= \Vert\tildex\Vert^2 \sqrt{\kappa(P)} \sqrt{\frac{1+4l^2}{\mu^2}}.
			\end{align*}
Further, we have when fixing $\beta/\alpha$, $\rho<1$ when $\alpha$ is small enough. 

\end{lemma}

Inexact gradient method with this type of gradient error has already been shown to have convergence guarantee, according to some related studies (e.g. \cite[Sec. IV-D]{qu2017harnessing}). For completeness we give a result of convergence for projected inexact gradient method in the following Lemma~\ref{lem:gd_inexact}, whose proof is deferred to Appendix-\ref{subsec:gd_inexact}.
\begin{lemma}\label{lem:gd_inexact}
	Recall that by Lemma~\ref{lem:L_property}, $\SSS$ is a concave and $\frac{\Vert\tildex\Vert^2}{\mu}$-smooth function and $\SSS$ is upper bounded. Consider the following algorithm
	$$ \vv{\lambda}(t+1) = [\vv{\lambda}(t) + \gamma \vv{g}(t)]^+$$
	where $[\cdot]^+$ is the projection onto the nonnegative orthant $\R_{\geq 0}^{2n}$, and $\vv{g}(t) $ is an inexact gradient that satisfies
	$$ \Vert \vv{g}(t) - \nabla \SSS(\vv{\lambda}(t))\Vert\leq C_1 \rho^{t} + C_2 \sum_{k=0}^{t-1}\rho^{t-1-k} \Vert \vv{\lambda}(k+1) -\vv{\lambda}(k)\Vert $$
	for $C_1,C_2>0$ and $\rho \in (0,1)$. Then, when $\gamma< \min(\frac{\mu}{2\Vert\tildex\Vert^2}, \frac{\Vert\tildex\Vert^2(1-\rho)^2}{2\mu C_2^2})$, we will have (i), $ \vv{\lambda}(t+1) - \vv{\lambda}(t) \rightarrow 0$, (ii) $\SSS (\vv{\lambda}(t))$ is lower bounded. 
\end{lemma}

The above lemma shows when $\gamma$ is small enough, $\vv{\lambda}(t+1) - \vv{\lambda}(t)\rightarrow 0$. By Lemma~\ref{lem:grad_error}, this further implies $\vv{\hatq}(t) - \vv{\hatq}^*(t) \rightarrow 0$ and $\vv{q}(t) - \vv{\hatq}^*(t)\rightarrow 0$. Since $S(\vv{\lambda}(t))$ is lower bounded, we have $\vv{\lambda}(t)$ is a bounded sequence by Lemma \ref{lem:L_property} (d). Then, sequence $(\vv{q}(t),\vv{\lambda}(t))$ is bounded and hence has a limit point $(\vv{q}_{lim} ,\vv{\lambda}_{lim})$, with subsequence $(\vv{q}(t_k),   \vv{\lambda}(t_k)   )$ converging to $(\vv{q}_{lim}, \vv{\lambda}_{lim} )$. Also since $\vv{\hatq}^*(t_k) = \vv{\hatq}^*(\vv{\lambda}(t_k))\rightarrow \vv{\hatq}^*(\vv{\lambda}_{lim})$, we have $\vv{q}_{lim} = \vv{\hatq}^*(\vv{\lambda}_{lim})$. We next show that $\vv{\lambda}_{lim}$ must be a maximizer of $\SSS(\vv{\lambda})$ over $\R^{2n}_{\geq 0}$. We have,
\begin{align*}
\lim_{k\rightarrow\infty} \vv{\lambda}(t_k+1) &=\lim_{k\rightarrow\infty} [ \vv{\lambda}(t_k) + \gamma(\tildex \vv{q}(t_k)  - \vv{v}_b) ]^+ \\
&= [\vv{\lambda}_{lim} +\gamma(\tildex \vv{q}_{lim}  - \vv{v}_b)  ]^+\\
&= [\vv{\lambda}_{lim} +\gamma(\tildex \vv{\hatq}^*(\vv{\lambda}_{lim})  - \vv{v}_b)  ]^+\\
&= [\vv{\lambda}_{lim}+\gamma \nabla \SSS(\vv{\lambda}_{lim})  ]^+.
\end{align*}  
Also since $\vv{\lambda}(t_k+1) - \vv{\lambda}(t_k) \rightarrow 0$, we have $\vv{\lambda}(t_k+1) \rightarrow\vv{\lambda}_{lim}$, and hence,
\begin{align}
\vv{\lambda}_{lim} = [\vv{\lambda}_{lim} +\gamma \nabla S(\vv{\lambda}_{lim})  ]^+.
\end{align}
This implies that, for any $\vv{\lambda}\in\R^{2n}_{\geq 0}$, 
\begin{align*}
&\nabla \SSS(\vv{\lambda}_{lim})^T(\vv{\lambda} - \vv{\lambda}_{lim} ) \\
&= \frac{1}{\gamma}(\vv{\lambda}_{lim}+ \gamma \nabla \SSS(\vv{\lambda}_{lim}) -  [ \vv{\lambda}_{lim}+ \gamma \nabla \SSS(\vv{\lambda}_{lim}) ]^+   )^T \\
&\quad \cdot (\vv{\lambda} - [ \vv{\lambda}_{lim}+ \gamma \nabla \SSS(\vv{\lambda}_{lim}) ]^+ )\\
&\leq 0
\end{align*}
where we have used the projection theorem. Since $\SSS$ is a concave function, we have $\vv{\lambda}_{lim}$ is indeed a maximizer of $\SSS$ over $\R^{2n}_{\geq 0}$. By Lemma~\ref{lem:L_property} (e), we have $\vv{q}_{lim} = \vv{\hatq}^*(\vv{\lambda}_{lim})$ is the unique solution of the original optimization problem (\ref{eq:opt}). This also shows the accumulation point of the bounded sequence $\vv{q}(t)$ is unique. So $\vv{q}(t)$ must converge to the unique solution of optimization problem (\ref{eq:opt}).

\begin{remark}\label{rem:stepsize}
	We here summarize what are the theoretic step size requirements for \algoname\ to converge. First we are given $l$, $\mu$, $c$, $Y$ as part of the problem parameters. Then, fix the ratio $\beta/\alpha$ to be any positive real number $\eta$, based on which the constants $a,\tau,\nu,C_2$ and matrix $P_0,P$ in Lemma~\ref{lem:grad_error} can be determined. Henceforth $\rho = e^{-\frac{\tau\alpha}{2}} + \frac{\alpha^2 \nu^2 \kappa(P_0) }{2}$ will only depend on $\alpha$. Notice that $\rho$ as a function of $\alpha$ equals $1$ when $\alpha = 0$, and have negative derivative at $\alpha = 0$, so to guarantee $\rho<1$ we only need to ensure $\alpha$ is small enough. Specifically, we can check that any $0<\alpha< \min(\frac{1}{\tau}, \frac{\tau}{2\nu^2\kappa(P_0)})$ ensures $\rho<1$. With $\alpha$ determined, and recall $\beta/\alpha$ is fixed to be $\eta$, we have $\beta = \eta\alpha$. Now $C_2$ and $\rho$ have both been determined, any $0<\gamma<\min(\frac{\mu}{2\Vert\tildex\Vert^2}, \frac{\Vert\tildex\Vert^2(1-\rho)^2}{2\mu C_2^2})$ will suffice. 
\end{remark}
\subsection{Proof of Lemma~\ref{lem:L_property}}	\label{appendix:proof_L}
	\noindent	\textbf{Proof of (a)}. We can show when fixing $\vv{\lambda}$, equation $\nabla_{\vv{\hatq}}\LL(\vv{\hatq},\vv{\xi},\vv{\lambda})=0,\nabla_{\vv{\xi}}\LL(\vv{\hatq},\vv{\xi},\vv{\lambda}) = 0$ is equivalent  to the KKT condition of the following optimization problem (where $\vv{\lambda}$ is fixed).
	\begin{subeqnarray}\label{eq:opt_fix_lambda}
		\min_{\vv{\hatq}} f(\vv{\hatq}) + \vv{\lambda}^T\tildex \vv{\hatq} \\ 
		s.t. \qquad  \unq_i\leq \hatq_i\leq\barq_i 
	\end{subeqnarray}
	The proof is almost identical to the proof of Lemma~\ref{lem:saddle}, and is hence omitted. Moreover, since problem~(\ref{eq:opt_fix_lambda}) has a unique primal-dual solution pair, we have the saddle point of $\LL(\cdot,\cdot,\vv{\lambda})$ is unique. Lastly, we have ${\hatq}^*_i(\vv{\lambda}) \in[\unq_i,\barq_i]$ since $\vv{\hatq}^*(\vv{\lambda})$ meets the constraint of (\ref{eq:opt_fix_lambda}).

	\noindent\textbf{Proof of (b).} 	It is easy to check the following relation
	$$\nabla_{\vv{\hatq}} K (\vv{\xi},\vv{\hatq})  = c \nabla_{\vv{\xi}} K (\vv{\xi},\vv{\hatq}) +\vv{\xi} .$$
	Then, since $\nabla_{\vv{\xi}}\LL(\vv{\hatq}^*(\vv{\lambda}),\vv{\xi}^*(\vv{\lambda}),\vv{\lambda}) = \nabla_{\vv{\xi}} K (\vv{\xi}^*(\vv{\lambda}),\vv{\hatq}^*(\vv{\lambda}))  = 0$, we have
	\begin{align*}
	0&=\nabla_{\vv{\hatq}} \LL(\vv{\hatq}^*(\vv{\lambda}),\vv{\xi}^*(\vv{\lambda}),\vv{\lambda} ) \\
	&=\nabla f(\vv{\hatq}^*(\vv{\lambda})) + \nabla_{\vv{\hatq}} K (\vv{\xi}^*(\vv{\lambda}),\vv{\hatq}^*(\vv{\lambda}))  + \tildex^T \vv{\lambda}\\
	&= \nabla f(\vv{\hatq}^*(\vv{\lambda})) + \vv{\xi}^*(\vv{\lambda}) + \tildex^T \vv{\lambda}.
	\end{align*}
	Writing down the same equation for $\vv{\lambda}'$, and taking the difference of the two, we have,
	{\small		\begin{align}
		&	\nabla f(\vv{\hatq}^*(\vv{\lambda})) - \nabla f(\vv{\hatq}^*(\vv{\lambda}')) +  \vv{\xi}^*(\vv{\lambda}) - \vv{\xi}^*(\vv{\lambda}')   + \tildex^T(\vv{\lambda} - \vv{\lambda}') = 0\label{eq:L:nabla_difference}
		\end{align}}Next, we claim that for each $i$, 
	\begin{align}
	(\xi_i^*(\vv{\lambda}) - \xi_i^*(\vv{\lambda}'))(\hatq_i^*(\vv{\lambda}) - \hatq_i^*(\vv{\lambda}')) \geq 0 \label{eq:L:innerproduct_positive}
	\end{align}
	This is because, if $\xi_i^*(\vv{\lambda})>0$, $\xi_i^*(\vv{\lambda}')>0$, then we must have $\hatq_i^*(\vv{\lambda}) = \hatq_i^*(\vv{\lambda}')  =\barq_i $ (cf. \eqref{eq:saddle:cs_1}), and hence \eqref{eq:L:innerproduct_positive} is true; if $\xi_i^*(\vv{\lambda})>0$, $\xi_i^*(\vv{\lambda}')\leq0$, then $\hatq_i^*(\vv{\lambda}) =\barq_i$, $ \hatq_i^*(\vv{\lambda}')  \leq\barq_i$ and \eqref{eq:L:innerproduct_positive} is still true. Other scenarios follow similarly. Using \eqref{eq:L:innerproduct_positive}, and take inner product between the left hand side of \eqref{eq:L:nabla_difference} and $\vv{\hatq}^*(\vv{\lambda}) - \vv{\hatq}^*(\vv{\lambda}') $, we have,
	\begin{align}
	0 &=  \langle \nabla f(\vv{\hatq}^*(\vv{\lambda})) - \nabla f(\vv{\hatq}^*(\vv{\lambda}')) , \vv{\hatq}^*(\vv{\lambda}) - \vv{\hatq}^*(\vv{\lambda}') \rangle \nonumber \\
	&\quad + \langle   \vv{\xi}^*(\vv{\lambda}) -\vv{\xi}^*(\vv{\lambda}'), \vv{\hatq}^*(\vv{\lambda}) - \vv{\hatq}^*(\vv{\lambda}')\rangle \nonumber \\
	&\quad + \langle \tildex^T(\vv{\lambda}-\vv{\lambda}'), \vv{\hatq}^*(\vv{\lambda}) - \vv{\hatq}^*(\vv{\lambda}')\rangle \label{eq:L:nabla_q_zero}\\
	&\geq \mu \Vert \vv{\hatq}^*(\vv{\lambda}) - \vv{\hatq}^*(\vv{\lambda}')  \Vert^2 - \Vert \tildex\Vert \Vert \vv{\lambda}-\vv{\lambda}'\Vert\Vert  \vv{\hatq}^*(\vv{\lambda}) - \vv{\hatq}^*(\vv{\lambda}')\Vert\nonumber
	\end{align}
	where we have used that $f$ is $\mu$-strongly convex. This implies 
	\begin{align*}
	\Vert \vv{\hatq}^*(\vv{\lambda}) - \hatq^*(\vv{\lambda}')  \Vert \leq \frac{\Vert\tildex\Vert }{\mu} \Vert \vv{\lambda} - \vv{\lambda}'\Vert .
	\end{align*}
	Further, by \eqref{eq:L:nabla_difference},
	{\small	\begin{align*}
		\Vert \vv{\xi}^*(\vv{\lambda}) - \vv{\xi}^*(\vv{\lambda}')\Vert &= \Vert \nabla f(\vv{\hatq}^*(\vv{\lambda})) - \nabla f(\vv{\hatq}^*(\vv{\lambda}'))  + \tildex^T(\vv{\lambda} - \vv{\lambda}') \Vert \\
		&\leq l\Vert \vv{\hatq}^*(\vv{\lambda}) - \vv{\hatq}^*(\vv{\lambda}')\Vert + \Vert\tildex\Vert \Vert \vv{\lambda}-\vv{\lambda}'\Vert\\
		&\leq 2\frac{l  }{\mu}\Vert\tildex\Vert  \Vert \vv{\lambda} - \vv{\lambda}'\Vert.
		\end{align*}}
	
	\noindent	\textbf{Proof of (c).} Clearly, $S(\vv{\lambda}) = \LL(\vv{\hatq}^*(\vv{\lambda}), \vv{\xi}^*(\vv{\lambda}),\vv{\lambda})$, and 
	$$\nabla S(\vv{\lambda}) = \nabla_{\vv{\lambda}}\LL(\vv{\hatq}^*(\vv{\lambda}), \vv{\xi}^*(\vv{\lambda}),\vv{\lambda} ) = \tildex \vv{\hatq}^*(\vv{\lambda}) - \vv{v}_b.$$
	In (\ref{eq:L:nabla_q_zero}),  we have shown the first two terms of the RHS of (\ref{eq:L:nabla_q_zero}) are nonnegative and hence $\langle \tildex^T(\vv{\lambda} - \vv{\lambda}') ,\vv{\hatq}^*(\vv{\lambda}) - \vv{\hatq}^*(\vv{\lambda}')\rangle\leq 0$. 
	Therefore,
	\begin{align*}
	\langle\nabla \SSS(\vv{\lambda}) - \nabla \SSS(\vv{\lambda}'),\vv{\lambda}-\vv{\lambda}'\rangle&= \langle \tildex (\vv{\hatq}^*(\vv{\lambda}) - \vv{\hatq}^*(\vv{\lambda}')),\vv{\lambda} - \vv{\lambda}'\rangle \\
	&= \langle \tildex^T(\vv{\lambda} -\vv{ \lambda}') ,\vv{\hatq}^*(\vv{\lambda}) - \vv{\hatq}^*(\vv{\lambda}')\rangle \\
	& \leq 0.
	\end{align*}
	Therefore, $\SSS$ is a concave function. To show the smoothness of $\SSS$, we have 
	\begin{align*}
	\Vert\nabla \SSS(\vv{\lambda}) - \nabla \SSS(\vv{\lambda}')\Vert &\leq \Vert\tildex\Vert \Vert \vv{\hatq}^*(\vv{\lambda}) - \vv{\hatq}^*(\vv{\lambda}')\Vert\\
	&\leq \frac{\Vert\tildex\Vert^2}{\mu} \Vert \vv{\lambda }- \vv{\lambda}'\Vert.
	\end{align*}
	
	\noindent		\textbf{Proof of (d).} Let $\vv{\hatq}^0\in\R^n$ be a feasible solution to \eqref{eq:opt} and meets \eqref{eq:vol_con}  with strict inequality. Then, for $\vv{\lambda} \in\R^{2n}_{\geq 0}$, noticing $\blambda_i(v_i(\vv{\hatq}^0) -\bv_i) \leq 0$ and $\ulambda_i(\uv_i - v_i(\vv{\hatq}^0))\leq 0$, we have,
	\begin{align*}
	&	\min_{\vv{\hatq}\in\R^n}	\LL(\vv{\hatq},\vv{\xi},\vv{\lambda})  \leq \LL(\vv{\hatq}^0,\vv{\xi},\vv{\lambda}) \\
	&= f(\vv{\hatq}^0) + \vv{\blambda}^T(\vv{v}(\vv{\hatq}^0) - \vv{\bar{v}})   + \vv{\ulambda}^T(\vv{\uv} - \vv{v}(\vv{\hatq}^0)) + K(\vv{\xi},\vv{\hatq}^0)\\
	&\leq f(\vv{\hatq}^0) + K(\vv{\xi},\vv{\hatq}^0).
	\end{align*}
	Since $\hatq^0_i \in [\unq_i,\barq_i]$, we have $K_i( 0 , \hatq^0_i) = 0$ and $\frac{\partial K_i(0,\hatq^0_i) }{\partial \xi_i} = 0$. This, together with the fact that $K_i(\xi_i, \hatq^0_i)$ is concave (in $\xi_i$), we have $K_i(\xi_i,\hatq^0_i)\leq 0, \forall \xi_i\in\R$. Therefore, $\min_{\vv{\hatq}\in\R^n}	\LL(\vv{\hatq},\vv{\xi},\vv{\lambda})  \leq f( \vv{\hatq}^0), \forall \vv{\xi}\in\R^n$ and hence, $\SSS(\vv{\lambda})  \leq f(\vv{\hatq}^0)$ (for all $\vv{\lambda} \in\R^{2n}_{\geq 0}$). Therefore $\SSS(\vv{\lambda})$ is bounded over $\R^{2n}_{\geq 0}$.
	
	Further, assuming $\vv{\lambda} \in Level_a $, then still using the fact $K(\vv{\xi},\vv{\hatq}^0)\leq 0$, we have
	\begin{align*}
	a&\leq \SSS(\vv{\lambda}) = \max_{\vv{\xi}\in\R^n} \min_{\vv{\hatq}\in\R^n}\LL(\vv{\hatq},\vv{\xi},\vv{\lambda})\\
	&\leq \max_{\vv{\xi}\in\R^n} \LL(\vv{\hatq}^0, \vv{\xi},\vv{\lambda})\\
	&\leq f(\vv{\hatq}^0) + \vv{\blambda}^T(\vv{v}(\vv{\hatq}^0) - \vv{\bar{v}})   + \vv{\ulambda}^T(\vv{\uv} -\vv{v}(\vv{\hatq}^0))
	\end{align*}
	Define $\epsilon = - \max_{i=1,\ldots,n} \max(v_i(\vv{\hatq}^0) - \bv_i, \uv_i - v_i(\vv{\hatq}^0)) $. We have $\epsilon>0$ since $\vv{\hatq}^0$ satisfies the voltage constraint with strict inequality. Hence, $a\leq f(\vv{\hatq}^0) -\epsilon \sum_{i=1}^n(\blambda_i+\ulambda_i)$ and 
	$$ \sum_{i=1}^n(\blambda_i+\ulambda_i) \leq \frac{1}{\epsilon}(f(\vv{\hatq}^0) - a).$$
	This together with the fact that $\vv{\lambda}\in\R^{2n}_{\geq 0}$ shows that the level set $Level_a$ is bounded. 
	
	\noindent		\textbf{Proof of (e). } By definition,  $(\vv{\hatq}^*(\vv{\lambda}^*), \vv{\xi}^*(\vv{\lambda}^*),\vv{\lambda}^*)$ is a saddle point of $\LL(\cdot,\cdot,\cdot)$. Then by Lemma~\ref{lem:saddle} we have $\vv{\hatq}^*(\vv{\lambda}^*)$ is the unique optimizer of problem~\eqref{eq:opt}. 

\subsection{Bounded Gradient Error, Proof of Lemma~\ref{lem:grad_error}}\label{subsec:bound_gd_error}

  Notice that by Lemma~\ref{lem:L_property} (a), $\vv{\hatq}^{*}(t)$ lies within the capacity constraint. Moreover, since $\vv{q}(t)$ is the projection of $\vv{\hatq}(t)$ onto the capacity constraint (cf. \eqref{eq:algorithm_q}), by projection theorem we have $\Vert \vv{q}(t) - \vv{\hatq}^*(t)\Vert\leq \Vert \vv{\hatq}(t) - \vv{\hatq}^*(t)\Vert $. Hence the gradient error satisfies
  \begin{align}
  \epsilon(t) \leq \Vert \tildex\Vert \Vert \vv{\hatq}(t) - \vv{\hatq}^*(t)\Vert. \label{eq:gd_error}
  \end{align}

By (\ref{eq:gd_error}), to bound the gradient error we only need to bound $\Vert \vv{\hatq}(t) - \vv{\hatq}^*(t)\Vert $. To this end, we now analyze what happens if we conduct one step of update for $\vv{\hatq}(t)$, $\vv{\xi}(t)$ (Eq. (\ref{eq:algorithm_hatq}) and (\ref{eq:algorithm_xi})). If we think of $\vv{\lambda}(t)$ as fixed, then the one step update for $\vv{\hatq}(t)$, $\vv{\xi}(t)$ is simply a primal-dual gradient update step for function $\LL(\cdot,\cdot,\vv{\lambda}(t))$.  
By our recent work on the geometric convergence rate of primal-dual gradient algorithm \cite{qu2018exponential}, after one step of update for $\vv{\hatq}(t)$, $\vv{\xi}(t)$ (Eq. (\ref{eq:algorithm_hatq}) and (\ref{eq:algorithm_xi})), the distance between $(\vv{\hatq}(t+1),  \vv{\xi}(t+1))$ and the unique saddle point of $\LL(\cdot,\cdot,\vv{\lambda}(t))$,  $(\vv{\hatq}^*(\vv{\lambda}(t)),\vv{\xi}^*(\vv{\lambda}(t))) = (\vv{\hatq}^*(t),\vv{\xi}^*(t))$ will shrink at least by a fixed ratio compared to the distance between $(\vv{\hatq}(t),\vv{\xi}(t))$ and $(\vv{\hatq}^*(t),\vv{\xi}^*(t))$. Here the distance is measured by a weighted norm. 
Formally, we stack $\vv{\hatq}(t), \vv{\xi}(t)$ into one large vector $\vv{z}(t) = [\vv{\hatq}(t)^T,\vv{\xi}(t)^T ]^T$, and similarly define $\vv{z}^*(t) = [\vv{\hatq}^*(t)^T, \vv{\xi}^*(t)^T]^T$. Then, we have the following lemma, which is a consequence of the results in \cite{qu2018exponential}. We will provide the detailed derivations in Appendix-\ref{subsec:primaldual}. 



\begin{lemma}\label{lem:primal_dual}
		Define $\mu'=\mu \underline{\sigma}(Y),l' = l\bar{\sigma}(Y)$. Further define constants,
	\begin{align*}
	a &= 20 l' [\max( \frac{c \bar{\sigma}(Y)}{\mu'}, \frac{l'}{\mu'} )]^2 [\max(\frac{\beta}{\alpha l' c},\frac{l'}{\mu'} ) ]^2 \kappa(Y)\\
	\tau& = \frac{\beta }{2\alpha } \frac{\underline{\sigma}(Y)}{a}\\
	\nu &= \max( \sqrt{2(l'+ \bar{\sigma}(Y)c)^2 + 2 \frac{\beta^2}{\alpha^2} \bar{\sigma}(Y)},\sqrt{2 \bar{\sigma}(Y) + 8 \frac{\beta^2}{\alpha^2 c^2}  }).
	\end{align*}
	Then we define matrix, $$P_0 = \left[\begin{array}{cc}
	\frac{\beta}{\alpha}a  I         & \frac{\beta}{\alpha} Y^{1/2}\\
	\frac{\beta}{\alpha} Y^{1/2} & a I
	\end{array}\right], P = \left[\begin{array}{cc}
	\frac{\beta}{\alpha} a Y^{-1}         & \frac{\beta}{\alpha} I\\
	\frac{\beta}{\alpha} I  & a I
	\end{array}\right]$$
	and constant,
	$$\rho = (e^{-\frac{\tau\alpha}{2}} + \frac{\alpha^2 \nu^2 \kappa(P_0) }{2} ).$$
Then, we have, 
	\begin{align}
		\Vert \vv{z}(t+1) - \vv{z}^*(t)\Vert_P \leq \rho \Vert \vv{z}(t) - \vv{z}^*(t)\Vert_P \label{eq:primal_dual_onestep}
	\end{align} 
where norm $\Vert \cdot \Vert_P$ is defined as  $\Vert \vv{z}\Vert_P = \sqrt{\vv{z}^TP\vv{z}}$. Further, fixing $\beta/\alpha$, we have $\rho<1$ when $\alpha$ is small enough.
\end{lemma} 

Lemma~\ref{lem:primal_dual}, specifically eq. (\ref{eq:primal_dual_onestep}), says that that $\vv{z}(t+1)$ gets closer to the equilibrium point $\vv{z}^*(t)$ by at least a fixed ratio $\rho$ compared to $\vv{z}(t)$, where the distance is measured in a specially constructed norm $\Vert \cdot\Vert_P$. From (\ref{eq:primal_dual_onestep}) we can easily get,
\begin{align}
	&\Vert \vv{z}(t) - \vv{z}^*(t)\Vert_P \nonumber\\
	&\leq  \rho \Vert \vv{z}(t-1) - \vv{z}^*(t-1)\Vert_P + \Vert \vv{z}^*(t)-\vv{z}^*(t-1)\Vert_P \nonumber \\
	& \leq \rho^t \Vert \vv{z}(0) - \vv{z}^*(0)\Vert_P + \sum_{k=0}^{t-1} \rho^{t-1-k}\Vert \vv{z}^*(k+1) -\vv{z}^*(k)\Vert_P  \label{eq:z_difference}
\end{align}
We now convert (\ref{eq:z_difference}) into a bound on $\Vert \vv{\hat{q}}(t) - \vv{\hatq}^*(t)\Vert$. 
First we have
\begin{align}
	\Vert \vv{z}(t) - \vv{z}^*(t)\Vert_P \geq   \sqrt{\underline{\sigma}(P)} \Vert \vv{\hatq}(t) - \vv{q}^*(t)\Vert.  \label{eq:P_2_convert_lb}
\end{align}

Recall that $\vv{z}^*(t) = [\vv{\hatq}^*(\vv{\lambda}(t))^T, \vv{\xi}^*(\vv{\lambda}(t))^T]^T$ is the unique saddle point of $\LL(\cdot,\cdot,\vv{\lambda}(t))$. 
Using Lemma~\ref{lem:L_property} (b) (also using $\Vert\tildex\Vert = \sqrt{2}\Vert X\Vert $),
\begin{align}
	&\Vert \vv{z}^*(t+1) - \vv{z}^*(t)\Vert_P^2  \nonumber \\
	&\leq  \bar{\sigma}(P) \Big[ \Vert \vv{\hatq}^*(t+1) - \vv{\hatq}^*(t)\Vert^2 +  \Vert \vv{\xi}^*(t+1) - \vv{\xi}^*(t)\Vert^2 \Big]\nonumber\\
	&\leq \bar{\sigma}(P)  \overbrace{\Big[ 2  (\frac{\Vert X\Vert}{\mu})^2 + 8  ( \frac{l}{\mu})^2 \Vert X\Vert^2 \Big]}^{:=C_0^2}  \Vert\vv{\lambda}(t+1) - \vv{\lambda}(t)\Vert^2. \label{eq:P_2_convert_ub}
\end{align} 

%
Using (\ref{eq:P_2_convert_lb}) and (\ref{eq:P_2_convert_ub}) in (\ref{eq:z_difference}), we have,

\begin{align}
	\epsilon(t) & \leq  \Vert \tildex\Vert \Vert \vv{\hatq}(t) - \vv{\hatq}^*(t)\Vert \nonumber\\
	&\leq  \Vert \tildex\Vert  \frac{1}{\sqrt{\underline{\sigma}(P)} }	\Vert \vv{z}(t) - \vv{z}^*(t)\Vert_P \nonumber \\
	& \leq \overbrace{  \Vert \tildex\Vert  \sqrt{\kappa(P)}  \Vert \vv{z}(0) - \vv{z}^*(0)\Vert}^{:=C_1}\rho^t \nonumber  \\
	&+ \underbrace{ \Vert \tildex\Vert  \sqrt{\kappa(P)} C_0  }_{:=C_2} \sum_{k=0}^{t-1} \rho^{t-1-k} \Vert \vv{\lambda}(k+1) - \vv{\lambda}(k) \Vert \label{eq:gd_error_full}
\end{align}
which concludes the proof of Lemma~\ref{lem:grad_error}.

	\subsection{Inexact Gradient Method, Proof of Lemma~\ref{lem:gd_inexact}}\label{subsec:gd_inexact}

\begin{proof}
	 Let $l'' = \frac{\Vert\tildex\Vert^2}{\mu}$. Then $\SSS$ is $l''$-smooth. By concavity and smoothness, 
{\small	\begin{align}
	&\SSS(\vv{\lambda}(t+1)) - \SSS(\vv{\lambda}(t))\nonumber \\
	&\geq  \langle \nabla \SSS(\vv{\lambda}(t)), \vv{\lambda}(t+1) - \vv{\lambda}(t)\rangle - \frac{l''}{2} \Vert \vv{\lambda}(t+1) - \vv{\lambda}(t)\Vert^2\nonumber\\
	&= \langle \vv{g}(t),\vv{\lambda}(t+1)- \vv{\lambda}(t)\rangle + \langle \nabla \SSS (\vv{\lambda}(t)) - \vv{g}(t),\vv{\lambda}(t+1) -\vv{ \lambda}(t)\rangle\nonumber\\
	&\quad - \frac{l''}{2} \Vert \vv{\lambda}(t+1) - \vv{\lambda}(t)\Vert^2.\label{eq:grad_error:lambda_1}
	\end{align}} By projection property, we have, 
	$$	\langle \vv{\lambda}(t) - \vv{\lambda}(t+1), \vv{\lambda}(t) +  \gamma \vv{g}(t) - \vv{\lambda}(t+1)\rangle\leq 0 $$
	which implies
	\begin{align}
	\langle \vv{g}(t), \vv{\lambda}(t+1) - \vv{\lambda}(t)\rangle \geq \frac{1}{\gamma} \Vert \vv{\lambda}(t+1)- \vv{\lambda}(t)\Vert^2. \label{eq:grad_error:projection}
	\end{align} 

 Also notice that 
	\begin{align*}
	& \langle \nabla \SSS (\vv{\lambda}(t)) - \vv{g}(t), \vv{\lambda}(t+1) - \vv{\lambda}(t) \rangle \\
	& \geq - \frac{l''}{2}\Vert\vv{\lambda}(t+1)- \vv{\lambda}(t)\Vert^2 - \frac{1}{2l''} \Vert \nabla \SSS(\vv{\lambda}(t)) - \vv{g}(t) \Vert^2.
	\end{align*}
	
	Plugging the above and (\ref{eq:grad_error:projection}) into (\ref{eq:grad_error:lambda_1}), we have
	\begin{align*}
	&\SSS (\vv{\lambda}(t+1))-\SSS (\vv{\lambda}(t)) \\
	&\geq ( \frac{1}{\gamma} -l'' )\Vert\vv{\lambda}(t+1) - \vv{\lambda}(t)\Vert^2 - \frac{1}{2l''}\Vert \nabla \SSS (\vv{\lambda}(t)) - \vv{g}(t)  \Vert^2
	\end{align*}
	
	Summing up the above from $t=0,\ldots, \tau$, we get,
	\begin{align}
	\SSS(\vv{\lambda}(\tau +1)) - \SSS(\vv{\lambda}(0)) &\geq  ( \frac{1}{\gamma} - l'' )\sum_{t=0}^\tau \Vert\vv{\lambda}(t+1)-\vv{\lambda}(t)\Vert^2 \nonumber \\
	&\quad - \frac{1}{2 l'' }\sum_{t=0}^\tau \Vert \nabla \SSS (\vv{\lambda}(t)) - \vv{g}(t)\Vert^2\label{eq:grad_error:lambda}
	\end{align}
	Now we bound $\sum_{t=0}^\tau \Vert \nabla \SSS (\vv{\lambda}(t)) - \vv{g}(t)\Vert^2$. By the condition in this lemma, we have
	\begin{align*}
	&\Vert \nabla \SSS(\vv{\lambda}(t)) - \vv{g}(t)\Vert  \\
	&\leq C_1 \rho^{t} + C_2 \sum_{k=0}^{t-1}\rho ^{t-1-k} \Vert \vv{\lambda}(k+1) - \vv{\lambda}(k)\Vert = \langle \vv{\chi}, \vv{\nu}_t\rangle 
	\end{align*}
	where we define vector $\vv{\chi} = [C_1,C_2 \Vert \vv{\lambda}(1)-\vv{\lambda}(0)\Vert, \ldots,C_2 \Vert\vv{\lambda}(\tau) - \vv{\lambda}(\tau-1) \Vert ]^T\in\R^{\tau+1}$, and vector $\vv{\nu}_t = [\rho^{t},\rho^{t-1}, \ldots,\rho,1,0,\ldots,0]^T\in\R^{\tau+1}$ for $0\leq t\leq \tau$. Then 
	$$\sum_{t=0}^\tau \Vert  \nabla \SSS(\vv{\lambda}(t)) - \vv{g}(t) \Vert^2 \leq \sum_{t=0}^\tau \vv{\chi}^T(\vv{\nu}_t\vv{\nu}_t^T)\vv{\chi}=\vv{\chi}^T V \vv{\chi} $$
	where $V=\sum_{t=0}^\tau \vv{\nu}_t \vv{\nu}_t^T\in\R^{(\tau+1)\times (\tau+1)}$. Obviously $V$ is symmetric and positive semi-definite, with each entry being nonnegative. Now, for $1\leq i\leq j\leq \tau+1$, the $(i,j)$'th entry of $V$ is $V_{ij}= \sum_{t=j-1}^\tau \rho^{t+1-i}\rho^{t+1-j} < \rho^{j-i} \frac{1}{1-\rho^2} < \frac{\rho^{j-i}}{1-\rho} $ .
	Therefore, fixing $i$, 
	$\sum_{j=i}^{\tau+1} |V_{ij}| < \sum_{j=i}^{\tau+1} \frac{\rho^{j-i}}{1-\rho} < \frac{1 }{(1-\rho)^2}  $, $\sum_{j=1}^{i-1}|V_{ij}| = \sum_{j=1}^{i-1} V_{ji}< \sum_{j=1}^{i-1} \rho^{i-j}\frac{1}{1-\rho} < \frac{1}{(1-\rho)^2}   $. So, $$\Vert V\Vert_1 \leq \sup_{i}\sum_{j=1}^{\tau+1}|V_{ij}|\leq \frac{2}{(1-\rho )^2}.$$
	This implies that
	\begin{align*}
	&\sum_{t=0}^\tau \Vert \nabla \SSS (\vv{\lambda}(t)) - \vv{g}(t)\Vert^2 \\
	&\leq \frac{2}{(1-\rho )^2}\Vert\vv{\chi}\Vert^2\\
	&\leq  \frac{2}{(1-\rho)^2 }C_1^2 +  \frac{2}{(1-\rho)^2}C_2^2 \sum_{t=0}^{\tau -1} \Vert \vv{\lambda}(t+1) - \vv{\lambda}(t)\Vert^2.
	\end{align*}
	Plugging this into (\ref{eq:grad_error:lambda}), we have
{\small	\begin{align*}
	&\SSS(\vv{\lambda}(\tau+1)) - \SSS(\vv{\lambda}(0))\\
	&\geq  ( \frac{1}{\gamma} - l'' - \frac{C_2^2}{l'' (1-\rho)^2} ) \sum_{t=0}^\tau \Vert\vv{\lambda}(t+1)-\vv{\lambda}(t)\Vert^2 - \frac{C_1^2}{l'' (1-\rho)^2}.
	\end{align*}}So if we choose $\gamma< \min(\frac{1}{2l''}, \frac{l''(1-\rho)^2}{2C_2^2})$, we have $\frac{1}{\gamma} - l'' - \frac{C_2^2}{l'' (1-\rho)^2} >0 $. Since $\SSS$ is upper bounded, we have  $\sum_{t=0}^\infty \Vert \vv{\lambda}(t+1) - \vv{\lambda}(t)\Vert^2 <\infty$, which implies $\Vert \vv{\lambda}(t+1) - \vv{\lambda}(t) \Vert\rightarrow 0$. The above inequality also implies that $\SSS(\vv{\lambda}(\tau+1))$ is lower bounded regardless of $\tau$.
\end{proof}

\subsection{Proof of Lemma~\ref{lem:primal_dual}}\label{subsec:primaldual}

In this section, since we only consider one step update, we drop the dependence on $t$ in the notations. Specifically, we write $\vv{\lambda}(t)$ as $\vv{\lambda}$, $\vv{\hatq}^*(t)$ and $\vv{\xi}^*(t)$ as $\vv{\hatq}^*$ and $\vv{\xi}^*$, $\vv{\hatq}(t)$ and $\vv{\xi}(t)$ as $\vv{\hatq}$ and $\vv{\xi}$,  $\vv{\hatq}(t+1)$ and $\vv{\xi}(t+1)$ as $\vv{\hatq}^+$ and $\vv{\xi}^+$, and at last, $\vv{z}^*(t)$, $\vv{z}(t)$ and $\vv{z}(t+1)$ as $\vv{z}^* $, $\vv{z}$ and $\vv{z}^+$. 

We now do the following change of variable from $\vv{\hatq}$ to $\vv{y}$, $ \vv{y} = Y^{-1/2} \vv{\hatq}$ ($\vv{y}^+$, $\vv{y}^*$ are defined similarly), while $\vv{\xi}$ stays unchanged. Correspondingly, vector $\vv{z}$ becomes 
$$\vv{w}= \left[\begin{array}{cc}
Y^{-1/2}        & 0\\
0&  I
\end{array}\right] \vv{z} =  \left[\begin{array}{c}
Y^{-1/2} \vv{\hatq}\\
\vv{\xi}
\end{array}\right] $$
and $\vv{w}^+$, $\vv{w}^*$ are defined correspondingly.
We then rewrite the Lagrangian in the new variables, $\tilde{\LL}(\vv{y},\vv{\xi}) = \LL(Y^{1/2} \vv{y},\vv{\xi},\vv{\lambda})$ (we drop the dependence of $\tilde{\LL}$ on $\vv{\lambda}$ since throughout this section, $\vv{\lambda}$ is fixed). Then,  the saddle point of $\tilde{\LL}$ is simply $\vv{w}^*$. Recall the update equation from $\vv{z}$ to $\vv{z}^+$ are
\begin{align*}
\vv{\hatq}^+ &= \vv{\hatq} - \alpha Y \nabla_{\vv{\hatq}} \LL(\vv{\hatq},\vv{\xi},\vv{\lambda})\\
\vv{\xi}^+ &=\vv{\xi} + \beta\nabla_{\vv{\xi}} \LL(\vv{\hatq},\vv{\xi},\vv{\lambda}).
\end{align*}
Now we rewrite the above update equation in variable $\vv{w}$ and get, 
\begin{align}
\vv{y}^+ &= \vv{y} - \alpha Y^{1/2} \nabla_{\vv{\hatq}}\LL(\vv{\hatq},\vv{\xi},\vv{\lambda}) \nonumber\\
&= \vv{y} - \alpha \nabla_{\vv{y}} \tilde{\LL}(\vv{y},\vv{\xi}) \label{eq:primal_dual:w_p}\\
\vv{\xi}^+&= \vv{\xi} + \beta \nabla_{\vv{\xi}} \LL(\vv{\hatq},\vv{\xi},\vv{\lambda})\nonumber\\
&= \vv{\xi} + \beta \nabla_{\vv{\xi}}\tilde{\LL}(\vv{y},\vv{\xi}). \label{eq:primal_dual:w_d}
\end{align}
Next, notice that
$\tilde{\LL} (\vv{y},\vv{\xi}) = f(Y^{1/2}\vv{y}) + \vv{\lambda}^T(\tildex Y^{1/2} \vv{y} - \vv{v}_b ) + K(\vv{\xi}, Y^{1/2}\vv{y})$ is precisely the augmented Lagrangian for the following optimization problem,
\begin{align}
&\min_{\vv{y}} f(Y^{1/2}\vv{y}) + \vv{\lambda}^T(\tildex Y^{1/2} \vv{y} - \vv{v}_b)  \label{eq:primal_dual:opt} \\
\text{s.t. }&\qquad  \vv{\unq}\leq Y^{1/2} \vv{y}\leq \vv{\barq}. \nonumber
\end{align}
The objective in \eqref{eq:primal_dual:opt} is $\mu'$-strongly convex and $l'$-smooth, where $\mu' = \mu \underline{\sigma}(Y) $, $l' = l \bar{\sigma}(Y)$. Then, by \cite[Thm. 3, Lem. 4, Lem. 5]{qu2018exponential}, we have the update~\eqref{eq:primal_dual:w_p} \eqref{eq:primal_dual:w_d} is a contraction, which is formally stated in the following Proposition.
\begin{proposition}(\cite[Thm. 3, Lem. 4, Lem. 5]{qu2018exponential})
	Given $\mu',l'$, $c$ and matrix $Y$, define constants,
	\begin{align*}
	a &= 20 l' [\max( \frac{c \bar{\sigma}(Y)}{\mu'}, \frac{l'}{\mu'} )]^2 [\max(\frac{\beta}{\alpha l' c},\frac{l'}{\mu'} ) ]^2 \kappa(Y)\\
	\tau& = \frac{\beta }{2\alpha } \frac{\underline{\sigma}(Y)}{a}\\
	\nu &= \max( \sqrt{2(l'+ \bar{\sigma}(Y)c)^2 + 2 \frac{\beta^2}{\alpha^2} \bar{\sigma}(Y)},\sqrt{2 \bar{\sigma}(Y) + 8 \frac{\beta^2}{\alpha^2 c^2}  }).
	\end{align*}
	Then we define matrix, $$P_0 = \left[\begin{array}{cc}
	\frac{\beta}{\alpha}a  I         & \frac{\beta}{\alpha} Y^{1/2}\\
	\frac{\beta}{\alpha} Y^{1/2} & a I
	\end{array}\right], $$
	and constant,
	$$\rho = (e^{-\frac{\tau\alpha}{2}} + \frac{\alpha^2 \nu^2 \kappa(P_0) }{2} ).$$
	 Then we have,
	$$\Vert \vv{w}^+ - \vv{w}^*\Vert_{P_0} \leq \rho \Vert \vv{w} - \vv{w}^*\Vert_{P_0}$$
	and further, fixing the ratio $\beta/\alpha$, we have $\rho<1$ if $\alpha$ is small enough.
\end{proposition}
Since $\vv{w}$,$\vv{w}^+$ and $\vv{w}^*$ are linear transforms of $\vv{z}$,$\vv{z}^+$ and $\vv{z}^*$, we have 
$$\Vert \vv{z}^+ - \vv{z}\Vert_P \leq \rho \Vert \vv{z} - \vv{z}^*\Vert_P $$
where 
\begin{align*}
P &=  \left[\begin{array}{cc}
Y^{-1/2}      & 0\\
0&  I
\end{array}\right] \left[\begin{array}{cc}
\frac{\beta}{\alpha}a  I         & \frac{\beta}{\alpha} Y^{1/2}\\
\frac{\beta}{\alpha} Y^{1/2} & a I
\end{array}\right]  \left[\begin{array}{cc}
Y^{-1/2}      & 0\\
0&  I
\end{array}\right] \\
&=\left[\begin{array}{cc}
\frac{\beta}{\alpha} a Y^{-1}         & \frac{\beta}{\alpha} I\\
\frac{\beta}{\alpha} I  & a I
\end{array}\right].\end{align*}
This concludes the proof of Lemma~\ref{lem:primal_dual}.
}{}
\end{document}